\documentclass{siamart220329}

\def\la{\lambda}

\def\e{{\epsilon}}

\usepackage[left]{lineno}
\usepackage{graphicx}
\usepackage{subfigure}
\usepackage{amsmath}
\usepackage{booktabs} 
\usepackage{multirow}
\usepackage{array}
\usepackage{slashed}
\usepackage{amsmath,bm}
\usepackage{appendix}
\usepackage{epsfig}
\usepackage{fancyhdr}
\usepackage{graphicx}
\usepackage{multirow}
\usepackage{float}
\usepackage{amsmath,amssymb}
\usepackage{graphicx}
\usepackage{algorithm}
\usepackage{tcolorbox}
\usepackage{algorithmic}
\usepackage{tabularx} 
\usepackage{ragged2e} 
\usepackage{booktabs} 
\usepackage{tablefootnote}

\newcolumntype{L}{>{\RaggedRight\hangafter=1\hangindent=0em}X}
\RequirePackage{algorithm,algorithmic}
\newcommand{\tabincell}[2]{\begin{tabular}{@{}#1@{}}#2\end{tabular}}

\newtheorem{remark}{Remark}[section]

\title{A block Lanczos method for large-scale quadratic minimization problems with orthogonality constraints}
\author{Bo Feng\thanks{School of Mathematics, China University of Mining and Technology, 221116, Jiangsu, P.R. China. E-mail: {\tt bofeng@cumt.edu.cn}.}
       \and Gang Wu\thanks{Corresponding author. School of Mathematics,
China University of Mining and Technology, Xuzhou, 221116, Jiangsu, P.R. China.
E-mail: {\tt gangwu@cumt.edu.cn}.
 This author is supported by the National Natural Science Foundation of China under grant 12271518, the Key Research and Development Project of Xuzhou Natural Science Foundation under grant KC22288, and the Open Project of Key Laboratory of Data Science and Intelligence Education of the Ministry of Education under grant DSIE202203.
}}

\begin{document}
\maketitle
\begin{abstract}
Quadratic minimization problems with orthogonality constraints (QMPO) play an important role in many applications of science and engineering.
However, some existing methods may suffer from low accuracy or heavy workload for large-scale QMPO. Krylov subspace methods are popular for large-scale optimization problems. In this work, we propose a block Lanczos method for solving the large-scale QMPO. In the proposed method, the original problem is projected into a small-sized one, and the Riemannian Trust-Region method is employed to solve the reduced QMPO. Convergence results on the optimal solution, the optimal objective function value, the multiplier and the KKT error are established. Moreover, we give the convergence speed of optimal solution, and show that if the block Lanczos process terminates, then an exact KKT solution is derived.
Numerical experiments illustrate the numerical behavior of the proposed algorithm, and demonstrate that it is more powerful than many state-of-the-art algorithms for large-scale quadratic minimization problems with orthogonality constraints.
\end{abstract}

\section{Introduction}
We are interested in solving the following large-scale quadratic minimization problems with orthogonality constraints (QMPO)
\begin{equation}\label{1}
\min_{U\in {\mathbb{O}}^{n\times \ell}}\left\{f(U):={\rm tr}(U^THU)+2~{\rm tr}(U^TG)\right\},\quad n>\ell>1,
\end{equation}
where $H=H^T\in \mathbb{R}^{n\times n}$ is symmetric, $G\in \mathbb{R}^{n\times \ell}$, $\mathbb{O}^{n\times \ell}:=\{U\in\mathbb{R}^{n\times\ell}~|~U^TU=I\}$.
This problem of \eqref{1} arises from many practical problems such as
orthogonal least squares regression (OLSR) \cite{OLSR1}, large graph clustering \cite{GCSED}, multidimensional similarity structure analysis \cite[Chapter 19]{Borg}, the Maxbet problem from canonical correlation analysis \cite{Geer,Liu}, multi-view subspace clustering \cite{zpei}, and so on.

Indeed, the famous \emph{unbalanced Procrustes problem} \cite{Chu,Edelman,Park,Gower,Viklands,Zl1,Zl2,Z.L}
\begin{equation}\label{1.1}
\min_{U\in {\mathbb{O}}^{n\times \ell}}\|AU-B\|_{\rm F},
\end{equation}
is a special case of QMPO, where $A \in \mathbb{R}^{
m\times n}$ and $B \in \mathbb{R}^{
m\times \ell}$. By the first-order optimality conditions for unbalanced Procrustes problem \eqref{1.1} (cf. \cite[Theorem 3.8]{Chu} and \cite[Theorem 3.1]{Park}),
we have following first-order optimality necessary conditions on QMPO \eqref{1}. There are some other types of optimality conditions on \eqref{1}, for more details, refer to \cite{Chu,Park,Zl2,Z.L}.
\begin{theorem}\cite[Theorem 3.8]{Chu} and \cite[Theorem 3.1]{Park}\label{ThmKKT}
If $U \in  \mathbb{O}^{n\times \ell}$
is a local minimizer of \eqref{1}, then there is a symmetric matrix $\Lambda\in \mathbb{R}^{\ell\times\ell}$
such that
\begin{equation}\label{th1}
HU+U\Lambda=-G.
\end{equation}
If $U$ is a global minimizer of \eqref{1}, then
\begin{equation}\label{glopt}
U^THU +\Lambda = - U^TG = - G^TU \succcurlyeq \bf O.
\end{equation}
\end{theorem}

Moreover, if $n=\ell$, then \eqref{1} reduces to the \emph{balanced Procrustes problem} \cite{Golub,Gower,Hur}
\begin{equation}\label{eqn1.5}
 \min_{U\in\mathbb{O}^{\ell \times \ell}}{\rm tr}(U^TG).
\end{equation}
In this case, we have a closed-form solution of \eqref{eqn1.5} by using the SVD decomposition of $G$ \cite{Golub,Gower}.
And if $\ell=1$, then \eqref{1} reduces to the classical {\it trust-region subproblem} \cite{3,0,Conn,Feng,10,W.W,5,28,7,4,Z.L1}.


Unfortunately, there is no closed-form solution for \eqref{1} generally. Some necessary or sufficient conditions for local and/or global minimizer of
\eqref{1} were established in \cite{Chu,Park,Zl2,Z.L}. Many iterative methods have been developed for the more general {\it optimization problems with orthogonal constraints}, which can be applied to solve \eqref{1} directly. For instance, Absil \emph{et al.} proposed a Riemannian Trust-Region (RTR) algorithm \cite{Abs,Absbook} for optimizing a smooth function on a Riemannian manifold. In \cite{Jiang}, Jiang and Dai proposed a framework for a constraint preserving update scheme for optimization on Stiefel manifold. In \cite{Wen}, Wen and Yin applied the Cayley transform to preserve the orthogonal constraints and develop curvilinear search algorithms with lower flops compared to those based on projections and geodesics. In \cite{Yuan2}, structured quasi-Newton methods were studied for optimization problems with orthogonality constraints. In \cite{Yuan4}, Gao {\it et al.} proposed a proximal linearized augmented Lagrangian algorithm for solving optimization problems with orthogonality constraints. A first-order framework was proposed in \cite{Yuan1} for optimization problems with orthogonal constraints.

The generalized power iteration (GPI) is one of the most popular methods for \eqref{1} \cite{Niee}. However, this method often suffers from the difficulty of slow convergence, and more detailed analysis is desired for the convergence of GPI. Recently, a novel eigenvalue-based approach was proposed in \cite{Z.L} to solve the unbalanced Procrustes problem \eqref{1.1}. This method also applies to the QMPO problem. It was proven that \eqref{1} can be equivalently transformed into an eigenvalue minimization whose solution can be computed by solving a related eigenvector-dependent nonlinear eigenvalue problem. However, one has to solve an $n$-by-$n$ (possibly dense) symmetric eigenproblem in each iteration of this algorithm, and the algorithm may converge very slowly if there is no subspace speeding up.




To the best of our knowledge, there are few specialized methods for solving {\it large-scale} QMPO \eqref{1}. Some existing methods may suffer from low accuracy or heavy workload for large-scale QMPO. Krylov subspace method is a powerful tool for solving large-scale optimization problems \cite{BF,10,5,4,Zl3,Zl4,Z.L1}. As far as we know, it seems that there is no (block) Krylov subspace method for the large-scale QMPO \eqref{1} till now. In this paper, we propose a block Krylov subspace method to solve \eqref{1}, in which the large-scale QMPO \eqref{1} is reduced into a small-sized one by using projection techniques. Furthermore, we establish the convergence results on the optimal solution, the optimal objective function value, the multiplier, as well as the KKT error. We give the convergence speed of optimal solution, and show that if the block Lanczos process terminates, then an exact KKT solution is derived, which satisfies the first order optimality in Theorem \ref{ThmKKT} and also the necessary condition \eqref{glopt} for a global minimizer. Numerical experiments on both synthetic and real-world data sets demonstrate that the proposed algorithm is superior to many state-of-the-art approaches for solving the large-scale QMPO \eqref{1}.

This paper is organized as follows. In Section 2, we propose a block Lanczos method for solving the large-scale QMPO.
The convergence of the proposed method is established in Section 3.
Numerical experiments are performed in Section 4 to show the numerical behavior of the new algorithm.
Some concluding remarks are given in Section 5.
Throughout this paper, we denote by $(\cdot)^T$ the transpose of a matrix or vector,
by $\mathcal{R}(E)$ the range space of a matrix $E$, and by $E\otimes F$ the Kronecker product of $E$ and $F$. In this paper, $E\succcurlyeq{ \bf O}~(E\succ{ \bf O})$ implies that $E$ is symmetric semi-positive definite (positive definite). Let $F=[\bm f_1,\ldots,\bm f_p]\in \mathbb{R}^{p\times q}$, then
$$
{\rm vec}(F)=(\bm f^T_1,\ldots,\bm f^T_p)^T\in \mathbb{R}^{pq}.
$$
Let $\bm 0$, ${ \bf O}$ and $ I$ be the zero
vector, zero matrix and identity matrix, respectively, whose orders are clear from the context.

%

\section{A block Lanczos method for solving the large-scale QMPO}
In this section, we propose a block Lanczos method to solve \eqref{1}. Let $G=V_1K$ be the economized QR decomposition of $G$, where $V_1\in \mathbb{O}^{n\times\ell}$. As $H$ is a symmetric matrix, we use the $k$-step block Lanczos process \cite{Golub,2,Y.S} to generate an orthonormal basis
$\mathbf{V}\!_{k }$ for the block Krylov subspace
$$
\mathcal{K}_k:=\mathcal{K}_k(H,V_1)=span\{V_1,HV_1,\ldots,H^{k-1}V_1 \}.
$$
Moreover, we have the following relation for this process \cite{Golub,2,Y.S}
\begin{equation}\label{eq1653}
H\mathbf{V}\!_{k } = \mathbf{V}\!_{k }T_{k }+V_{k+1}N_{k}(E^{(k\ell)}_{\ell})^T,
\end{equation}
where  $\mathbf{V}\!_{k}:=[V_1,~V_2,\ldots,~V_k]\in \mathbb{R}^{n\times k\ell}$, $\mathbf{V}_{k }^T\mathbf{V}\!_{k }=I_{k\ell}$, $\mathbf{V}_{k }^TV_{k+1}={ \bf O}$, $V_{k+1}^TV_{k+1}=I_{\ell}$, $N_k\in \mathbb{R}^{\ell\times \ell}$ is upper triangular, and $E^{k\ell}_{\ell}$ denotes the last $\ell$ columns of the identity matrix $I_{k\ell}$. Here
\begin{equation}\label{eq1121}
T_{k }=\mathbf{V}_{k }^TH\mathbf{V}\!_{k }=\begin{pmatrix}
M_1  & N^T_1  &   &  & \\
N_1  & M_2  & N^T_2 & & \\
     &   \ddots  & \ddots &  \ddots \\
     &    & \ddots &     M_{k-1} & N^T_{k-1}\\
     &        &    &        N_{k-1}& M_k
\end{pmatrix}\in \mathbb{R}^{k\ell\times k\ell},
\end{equation}
is block tridiagonal,
with $M_j \in \mathbb{R}^{\ell\times \ell}$, and $N_j\in \mathbb{R}^{\ell\times \ell}$ being upper triangular, $j=1,2,\ldots,k$.




In the proposed method, \eqref{1} reduces to the following small-sized {\it constrained} problem:
\begin{equation}\label{2}
\min_{\substack{U\in \mathbb{O}^{n\times \ell}\\ \mathcal{R}(U)\subseteq \mathcal{K}_k}}\left\{f(U)={\rm tr}(U^THU)+2~{\rm tr}(U^TG)\right\}.
\end{equation}
Indeed, \eqref{2} can be equivalently rewritten as the following \emph{reduced} QMPO:
 \begin{equation}\label{3}
\min_{P\in \mathbb{O}^{k\ell\times \ell}}\left\{\widetilde{f}(P)={\rm tr}(P^TT_kP)+2~{\rm tr}(P^T{G}_k)\right\},
\end{equation}
where
 $$
 {G}_k= \mathbf{V}_k^TG\!=\!\mathbf{V}_k^TV_1K\!=\!
         \begin{pmatrix}K\\{\bf O}  \end{pmatrix}\in \mathbb{R}^{ k\ell\times\ell}.
 $$

Let $P_k$ be a solution to \eqref{3},
then
we use
\begin{equation}\label{eq2121}
U_k= \mathbf{V}\!_kP_k=\arg\min_{\substack{U\in \mathbb{O}^{n\times \ell}
\\ \mathcal{R}(U)\subseteq \mathcal{K}_k}}f(U)
\end{equation}
as an approximation to the optimal solution $U_*$, and
\begin{equation}\label{eqn2.5}
 f(U_k)=f(\mathbf{V}\!_kP_k)=\widetilde{f}(P_k)
\end{equation}
is an approximation to the optimal value $f(U_*)$.
{
By Theorem \ref{ThmKKT}, there is a symmetric matrix $\Lambda_k\in\mathbb{R}^{k\times k}$ such that
\begin{equation}\label{eq10}
T_kP_k+P_k\Lambda_k+G_k=0.
\end{equation}
}

Consequently, we reduce the large-scale QMPO
\eqref{1} to a $k\ell$-by-$k\ell$ small-sized one. In practice, one can exploit the Riemannian Trust-Region
(RTR) method \cite{Abs,Absbook} to solve \eqref{3}.
The proposed algorithm is given in Algorithm \ref{alg}. One refers to Section \ref{Sec4} for more details on practical implementations.
\renewcommand{\thealgorithm}{1}
\begin{algorithm}
\caption{A block Lanczos method for large-scale QMPO}
\label{alg}
\begin{algorithmic}[1]
\REQUIRE 
    $H\in \mathbb{R}^{n\times n}$, $G\in\mathbb{R}^{n\times \ell}$, and $k_{\max}$.
\ENSURE 
    $U_k$.\\
   \STATE Set $V_0={\bf O}\in \mathbb{R}^{n\times \ell}$,
    $N_0 = {\bf O}\in \mathbb{R}^{\ell\times \ell}$ and $k=0$;
    \STATE   Compute the economized QR decomposition: $G =V_{1}K$, where $V_1\in\mathbb{R}^{n\times\ell}$;
  \STATE  Let $M_1=V_1^THV_1$;
  \STATE   ${\bf whlie}$ $k\leq k_{\max}$
  \STATE    ~~~$k=k+1$;
  \STATE    ~~~Let $L_k=HV_k - V_k M_k - V_{k-1}N_{k-1}^T$;
  \STATE    ~~~Compute the economized QR decomposition: $L_k =V_{k+1}N_k$;
  \STATE    ~~~Let $M_{k+1}=V_{k+1}^THV_{k+1}$;~~~~
  \STATE    ~~~Solve the reduced QMPO
 \begin{equation*}
P_k =\arg\min_{P\in \mathbb{O}^{k\ell \times \ell}}\Big\{{\rm tr}(P^TT_k P)+2~{\rm tr}(P^T{G}_k)\Big\},
\end{equation*}
            ~~~where $T_k\in\mathbb{R}^{k\ell\times k\ell}$ is defined in \eqref{eq1121} and {$ {G}_k=
         \begin{pmatrix}K\\ {\bf O} \end{pmatrix}\in \mathbb{R}^{ k\ell\times\ell}$};
  \STATE    ~~{\bf if} the convergence criterion is satisfied~~ \%~~{\tt Refer to \eqref{eq1234}}
   \STATE    ~~~~$U_k= \mathbf{V}\!_kP_k$;
  \STATE    ~~{\bf end~if}
  \STATE    ${\bf end~whlie}$
\end{algorithmic}
\end{algorithm}

\section{Convergence analysis}
In this section, we show the convergence of Algorithm \ref{alg}. We first need the following three lemmas.
The first lemma follows from the definition of the Kronecker product and \cite[Theorem 4.4.5]{H.B}.
\begin{lemma}\label{lem1627}
Let $X\in \mathbb{R}^{s\times s}$ and $Y\in\mathbb{R}^{m\times m}$, then
$(I_m\!\otimes\! X)\!+\!(Y\!\otimes\! I_s)$  is nonsingular if and only if $X\!+\!\la_Y I$ is nonsingular, where $\la_Y$ is an eigenvalue of $Y$.
Moreover, if $X=X^T$ and $Y=Y^T$, then
\begin{subequations}
\begin{align}
&(I_m\otimes X)+(Y\otimes I_s)~{\rm is ~symmetric},\label{eq2031}\\
&\la_{\max}\big((I_m\otimes X)+(Y\otimes I_s)\big)=\la_{\max}(X) +\la_{\max}(Y),\label{eq2031a}\\
&\la_{\min}\big((I_m\otimes X)+(Y\otimes I_s)\big)=\la_{\min}(X) +\la_{\min}(Y).\label{eq2031b}
\end{align}
\end{subequations}
\end{lemma}

The second lemma is from \cite[Section 4.2, Problem 25]{H.B}.
\begin{lemma}\cite{H.B}\label{eq1614}
Let $C\in \mathbb{R}^{t\times s}$, $E\in \mathbb{R}^{p\times q}$,
 $X\in\mathbb{R}^{q\times s}$ and $Y\in\mathbb{R}^{p\times t}$. Then
\begin{equation}\label{eq1148}
{\rm tr}(C^TY^TEX)={\underline{\bm y}}^T(C\otimes E)\underline{\bm x}
\end{equation}
where
$\underline{\bm y}={\rm vec}(Y)$
and
$
\underline{\bm x}={\rm vec}(X)
$.
\end{lemma}

The third lemma is the polar decomposition of a full column rank matrix.
\begin{lemma} \cite[Theorem 8.1]{J.H}\label{le1435}
Let $Y\in \mathbb{R}^{p\times s}$ ($p\geq s$) with $rank(Y)=s$.
There exists a unique matrix $Q\in \mathbb{R}^{p\times s}$ with orthonormal columns and a unique symmetric positive definite matrix $S$ such that $Y= QS$.
 The matrix $S$ is given by $S=(Y^TY)^{\frac{1}{2}}$.
\end{lemma}

We are ready to consider the convergence of the proposed method.
Let $U_*$ be a global minimizer of \eqref{1}, then $U_*$ is also a local solution. It follows from Theorem \ref{ThmKKT} that, there is a symmetric matrix $\Lambda_* \in \mathbb{R}^{\ell \times \ell}$, such that \eqref{th1} holds.
Let the eigendecompositions of $H$ and $\Lambda_*$ be
\begin{subequations}
\begin{align}
H&=W D W^T=W ~diag(\mu_1,\mu_2,\ldots,\mu_{n})~W^T,~~{\rm with}~~\mu_1\geq\mu_2\geq\cdots\geq\mu_{n},\label{eq2143}\\
\Lambda_*&=Z\Gamma Z^T=Z ~diag(\gamma_1,\gamma_2,\ldots,\gamma_{\ell})~Z^T,~~{\rm with}~~\gamma_1\geq\gamma_2\geq\cdots\geq\gamma_{\ell},\label{eq2143b}
\end{align}
\end{subequations}
where $W\in\mathbb{R}^{n\times n}$ and $Z\!=\![\bm z_1,\bm z_2,\ldots,\bm z_{\ell}]\!\in\! \mathbb{R}^{\ell\times \ell}$ are orthonormal matrices.

In this paper,
we make the following assumption
\begin{equation}\label{asp1}
Assumption:~
\mathcal{H}_*:=(I_{\ell}\otimes H)+(\Lambda_*\otimes I_n)~~\rm is ~~nonsingular.
\end{equation}
Hence, it follows from Lemma \ref{lem1627} that $H+\gamma_i I$ are nonsingular, $i=1,2,\ldots,\ell$.
Consider the two index sets
\begin{align*}
\mathcal{I}&=\{i~|~H+\gamma_i I\succ{\bf O},~{\rm i.e.,}~\mu_n+\gamma_i>0,~{\rm where}~1\leq i\leq \ell\},\\
\mathcal{J}&=\{i~|~H+\gamma_i I~~{\rm is~ nonsingular~and~}\mu_n+\gamma_i<0,~{\rm where}~1\leq i\leq \ell\}.
\end{align*}
Thus, we have that $\mathcal{I}\cup \mathcal{J}=\{1,2,\ldots,\ell\}$ and $\mathcal{I}\cap\mathcal{J}=\varnothing$.
It was shown in \cite[Theorem 2.4]{Z.L} that $\mu_1+\gamma_n\geq0$. In other words, $H+\gamma_i I$ will never be a negative definite matrix, $i=1,2,\ldots,\ell$. Therefore, if $\mathcal{J}\neq \varnothing$, there is an integer $1\leq s_i\leq n$, such that
$$
\mu_n+\gamma_i\leq \cdots \leq
 \mu_{s_i}+\gamma_i<0<\mu_{s_i+1}+\gamma_i\leq \cdots \leq
 \mu_{1}+\gamma_i,~~{\rm for}~~i\in \mathcal{J}.
$$

Consider 
$$
\phi_i:=\frac{a_ib_i}{|(\mu_{s_i}+\gamma_i)(\mu_{s_i+1}+\gamma_i)|}\geq 1
~~{\rm for}~~i\in\mathcal{J},
$$
where
\begin{align*}
a_i&=\max\big\{-(\mu_n+\gamma_i),\mu_{1}-\mu_{s_i+1}-(\mu_{s_i}+\gamma_i) \big\},\\
b_i&=\max\big\{\mu_{1}+\gamma_i,\mu_{s_1+1}-\mu_{n}+\mu_{s_i}+\gamma_i \big\}.
\end{align*}
The definitions of $a_i$ and $b_i$ are due to the embedding of {$[\mu_n\!+\!\gamma_i,
 \mu_{s_i}\!+\!\gamma_i]\cup
 [\mu_{s_i+1}\!+\!\gamma_i,\mu_{1}\!+\!\gamma_i]\!\subseteq
 \![-a_i,\mu_{s_i}\!+\!\gamma_i]\cup[\mu_{s_i+1}\!+\!\gamma_i,b_i]$}
into intervals of equal lengths \cite[section 3.1]{Herzog}.

$\bullet$ First, we consider the distance between the optimal solution $U_*$ and the search subspace $\mathcal{K}_k$. Indeed,
a necessary condition for the convergence of the proposed method is that the distance tends to zero.

\begin{theorem}\label{thm1540}
Let
$$
\e_k = \min_{X\in \mathbb{R}^{k\ell\times \ell}}\|U_*-\mathbf{V}\!_kX\|_{\rm F}=\|(I-\mathbf{V}\!_k\mathbf{V}_k^T) U_*\|_{\rm F}.
$$
Under the previous notation and assumption, we have
\begin{equation}\label{eq1009}
\e_k\leq2\cdot\sqrt{\sum_{i\in \mathcal{I}} \left(\frac{\sqrt{\kappa_i}-1}{\sqrt{\kappa_i}+1}\right)^{2(k+1)}+\sum_{i\in \mathcal{J}} \left(\frac{\sqrt{\phi_i}-1}{\sqrt{\phi_i}+1} \right)^{k-1}},
\end{equation}
where $\kappa_i$ is the 2-condition number of $H+\gamma_i I$ for $i\in \mathcal{I}$.
\end{theorem}
\begin{proof}
From $HU_*+U_*\Lambda_*=-G$ and $\Lambda_*=Z\Gamma Z^T$, we have
\begin{equation*}
(H+\gamma_i I)U_*\bm z_i=-G\bm z_i~~{\rm for}~~i=1,2,\ldots,\ell.
\end{equation*}
Denote by ${\mathcal{K}}^{(i)}_k\!:=\!\mathcal{K}_k(H\!+\!\gamma_i I,G\bm z_i)$ and let $\mathcal{P}_{t}$ be the set of polynomials with degree no higher than $t$. It holds that
\begin{align*}
\min_{\bm y\in \mathcal{K}^{(i)}_k}\!\|U_*\bm z_i\!-\!\bm y\|_2
&=\min_{p\in \mathcal{P}_k}\|U_*\bm z_i\!-\!p(H\!+\!\gamma_i I)G\bm z_i\|_2\\
&=\min_{p\in \mathcal{P}_k}\|[I\!+\!p(H\!+\!\gamma_i I)(H\!+\!\gamma_i I)]U_*\bm z_i\|_2\\
&=\min_{\substack{h\in \mathcal{P}_{k+1}\\ h(0)=1}}\|h(H+\gamma_i I)U_*\bm z_i\|_2\leq\min_{\substack{h\in \mathcal{P}_{k+1}\\ h(0)=1}}\|h(H+\gamma_i I)\|_2\\
&\overset{\eqref{eq2143}}{=}\min_{\substack{h\in \mathcal{P}_{k+1}\\ h(0)=1}}\|h(D+\gamma_i I)\|_2.
\end{align*}

From the Assumption \eqref{asp1}, we have that $H+\gamma_i I$ are nonsingular, $i=1,2,\ldots,\ell$.
Then it follows from \cite[Section 3.1]{8} that, on one hand, if $i\in \mathcal{I}$,
\begin{align}\label{eq2036}
\min_{\bm y\in \mathcal{K}^{(i)}_k}\!\|U_*\bm z_i\!-\!\bm y\|_2&\leq\min_{\substack{h\in \mathcal{P}_{k+1}\\ h(0)=1}}\|h(D+\gamma_i I)\|_2
\nonumber\\
     &\leq\underset{\substack{h\in \mathcal{P}_{k+1}\\ h(0)=1}}{\min}~\underset{t\in [\mu_n\!+\gamma_i,\mu_1\!+\gamma_i]}{\max}|h(t)|\leq
2\left(\frac{\sqrt{\kappa_i}-1}{\sqrt{\kappa_i}+1} \right)^{k+1}.
\end{align}
On the other hand, if $i\in \mathcal{J}$,
\begin{align}
\min_{\bm y\in \mathcal{K}^{(i)}_k}\!\|U_*\bm z_i&-\bm y\|_2
\leq\min_{\substack{h\in \mathcal{P}_{k+1}\\ h(0)=1}}\|h(D+\gamma_i I)\|\leq\underset{\substack{h\in \mathcal{P}_{k+1}\\ h(0)=1}}{\min}\underset{t\in\mathcal{L}}{\max}|h(t)|\nonumber\\
\leq&
2\left(\frac{\sqrt{|a_ib_i|}\!-\!\sqrt{|(\mu_{s_i}+\gamma_i)(\mu_{s_i+1}+\gamma_i)|}
}{\sqrt{|a_ib_i|}+\sqrt{|(\mu_{s_i}+\gamma_i)(\mu_{s_i+1}+\gamma_i)|}} \right)^{\lfloor\frac{k+1}{2}\rfloor}
= 2\left(\frac{\sqrt{\phi_i}-1}{\sqrt{\phi_i}+1} \right)^{\lfloor\frac{k+1}{2}\rfloor},
\label{eq2100}
\end{align}
where $\mathcal{L}=[-a_i,\mu_{s_i}+\gamma_i]
\cup[\mu_{s_i+1}+\gamma_i,b_i]$ and
$\lfloor\cdot\rfloor$ stands for the integer part of a number.

From $G\bm z_i\in\mathcal{R}(V_1)$, we obtain
\begin{equation*}
{\mathcal{K}}^{(i)}_k=\mathcal{K}_k(H,G\bm z_i)
\!\subseteq\!\mathcal{K}_k(H,V_1)=\mathcal{K}_k~~{\rm for}~~i=1,2,\ldots,\ell.
\end{equation*}
So we have
 \begin{align}
\e^2_k=&\left\|(I\!-\!\mathbf{V}\!_k\mathbf{V}_k^T)U_*\right\|^2_{\rm F}
\overset{\eqref{eq2143b}}{=}\left\|(I\!-\!\mathbf{V}\!_k
\mathbf{V}_k^T)U_*Z\right\|^2_{\rm F}\nonumber\\
=&\sum_{i=1}^{\ell}\!\|(I\!-\!\mathbf{V}\!_k\mathbf{V}_k^T)U_*\bm z_i\|^2_2
=\sum_{i=1}^{\ell}\min_{\bm y\in \mathcal{K}_k}\!\|U_*\bm z_i\!-\!\bm y\|^2_2\nonumber\\
=&\sum_{i\in\mathcal{I}}\min_{\bm y\in \mathcal{K}_k}\!\|U_*\bm z_i\!-\!\bm y\|^2_2+\sum_{i\in\mathcal{J}}\min_{\bm y\in \mathcal{K}_k}\!\|U_*\bm z_i\!-\!\bm y\|^2_2\nonumber\\
\leq&\sum_{i\in\mathcal{I}}\min_{\bm y\in \mathcal{K}^{(i)}_k}\!\|U_*\bm z_i\!-\!\bm y\|^2_2+\sum_{i\in\mathcal{J}}\min_{\bm y\in \mathcal{K}^{(i)}_k}\!\|U_*\bm z_i\!-\!\bm y\|^2_2,\nonumber
\end{align}
which, together with \eqref{eq2036} and \eqref{eq2100}, yields \eqref{eq1009}.
\end{proof}
\begin{remark}
Theorem \ref{thm1540} shows that the rate at which $\e_k$ converges to 0
strictly relies on the distribution of the spectrum of $H+\gamma_i I, i=1,2,\ldots,\ell$.
In particular, the convergence rate of $\e_k$ is comparable to that of
conjugate gradient method provided that $(I_{\ell}\otimes H)+(\Lambda_*\otimes I_n)\succ\bf O$.
\end{remark}

$\bullet$ Second, we show the convergence of $f(U_k)$. To this aim, we consider the upper bound of $f(U_k)\!-\!f(U_*)$.

\begin{theorem}
Suppose that
$\|(I-\mathbf{V}\!_k\mathbf{V}_k^T)U_*\|_2<1$. Then
\begin{equation}\label{eq1724}
0\leq f(U_k)-f(U_*)\leq 2(\mu_1+\gamma_1)\cdot\e_k^2.
\end{equation}
\end{theorem}
\begin{proof}
For any $U\in\mathbb{O}^{n\times\ell}$, we show that
\begin{equation}\label{eq21.54}
0\leq f(U)-f(U_*)={\rm tr}[(U_*-U)^TH(U_*-U)]+{\rm tr}[(U_*-U)\Lambda_* (U_*-U)^T].
\end{equation}
Indeed,
\begin{align}
0\!&\leq\! f(U)\!-\!f(U_*)
={\rm tr}(U^THU)-{\rm tr}({U}_*^TH{U}_*)
+2~{\rm tr}[(U-{U}_*)^TG]\nonumber\\
&\overset{\eqref{th1}}{=}{\rm tr}(U^T\!HU)\!-\!{\rm tr}({U}_*^T\!H{U}_*)
+{\rm tr}[(U\!-\!{U}_*)^TG]\!-\!{\rm tr}[(U\!-\!U_*)^T(HU_*\!+\!U_*\Lambda_*)]\nonumber\\
&={\rm tr}(U^THU)\!+\!{\rm tr}[(U-U_*)^TG]\!-\!{\rm tr}(U^THU_*)+{\rm tr}[(U_*-U)^TU_*\Lambda_*]\nonumber\\
&={\rm tr}[(U-U_*)^T(HU+G)]+{\rm tr}[(U_*-U)^TU_*\Lambda_*]\nonumber\\
&={\rm tr}[(U-U_*)^TH(U-U_*)]\!+\!{\rm tr}[(U-U_*)^T(HU_*\!+\!G)]
\!+\!{\rm tr}[(U_*\!-\!U)^TU_*\Lambda_*]\nonumber\\
&\overset{\eqref{th1}}{=}{\rm tr}[(U-U_*)^TH(U-U_*)]+2~{\rm tr}[(U_*-U)^TU_*\Lambda_*]\nonumber\\
&={\rm tr}[(U\!-\!U_*)^T\!H(U\!-\!U_*)]\!+\!{\rm tr}[(U_*\!-\!U)^T\!(U_*\!-\!U)\Lambda_*]
  \!+\!{\rm tr}[(U_*\!-\!U)^T\!(U_*\!+\!U)\Lambda_*]\nonumber\\
&={\rm tr}[(U-U_*)^TH(U-U_*)]+{\rm tr}[(U_*-U)\Lambda_* (U_*-U)^T],\nonumber
\end{align}
where the last inequality is from the facts that $U\in \mathbb{O}^{n\times\ell}$, $\Lambda_*=\Lambda_*^T$, and
\begin{align*}
{\rm tr}[(U_*-U)^T(U_*+U)\Lambda]&={\rm tr}[(U_*^TU-U^T U_*)\Lambda_*]\\
&={\rm tr}(U_*^TU\Lambda_*)-{\rm tr}(U^T U_*\Lambda_*)\\
&={\rm tr}(U\Lambda_*U_*^T)-{\rm tr}( U_*\Lambda_*U^T)\\
&={\rm tr}(U\Lambda_*U_*^T)-{\rm tr}( U\Lambda_*U_*^T)=0,
\end{align*}
so we have \eqref{eq21.54}.

We are ready to prove \eqref{eq1724}.
It follows from \cite[Theorem 3.3.16 (c)]{H.B} that
\begin{align*}
\big|1-\sigma_{\ell}(\mathbf{V}\!_k\mathbf{V}^T_k\!U_*)\big|
=&\big|\sigma_{\ell}(U_*)
-\sigma_{\ell}(\mathbf{V}\!_k\mathbf{V}^T_kU_*)\big|\\
\leq&
\left\|U_*-\mathbf{V}\!_k\mathbf{V}^T_kU_*\right\|_2=\|(I\!-\!\mathbf{V}\!_k\mathbf{V}^T_k)U_*\|_2\!<\!1.
\end{align*}
Hence, $\sigma_{\ell}(\mathbf{V}\!_k\mathbf{V}^T_kU_*)>0$, i.e., $\mathbf{V}\!_k\mathbf{V}^T_kU_*\in \mathbb{R}^{n\times \ell}$ has full column rank. By Lemma \ref{le1435},
there is a unique orthonormal matrix $\widetilde{U}_k\in \mathbb{C}^{n\times \ell}$ and a  symmetric positive definite matrix $M$, such that
$
\mathbf{V}\!_k\mathbf{V}^T_kU_*\!\!=\!\widetilde{U}_k M
$.
Thus,  $ \mathcal{R}(\widetilde{U}_k)\subseteq
\mathcal{R}(\mathbf{V}_k)=\mathcal{K}_k$.
By \eqref{eq2121},
\begin{align}\label{eq21.57}
0&\leq f(U_k)-f(U_*)\leq f(\widetilde{U}_k)-f(U_*)\nonumber\\
             &\overset{\eqref{eq21.54}}{=}{\rm tr}[(\widetilde{U}_k-U_*)^TH(\widetilde{U}_k-U_*)]
              +{\rm tr}[(U_*-\widetilde{U}_k)
\Lambda_* (U_*-\widetilde{U}_k)^T]\nonumber\\
            &={\rm tr}[(\widetilde{U}_k-U_*)^T H(\widetilde{U}_k-U_*)]
                                    +{\rm tr}[\Lambda_* (U_*-\widetilde{U}_k)^T(U_*-\widetilde{U}_k)]\nonumber\\
            &\overset{\eqref{eq1148}}{=}\big[{\rm vec}(U_*-\widetilde{U}_k)\big]^T\cdot\mathcal{H}_*\cdot{\rm vec}\big(U_*-\widetilde{U}_k)\nonumber\\
            &\overset{\eqref{eq2031a}}{\leq} (\mu_1+\gamma_1)\cdot\left\|{\rm vec}(U_*-\widetilde{U}_k)\right\|_2^2
            =(\mu_1+\gamma_1)\cdot\left\|U_*-\widetilde{U}_k\right\|_{\rm F}^2.
\end{align}

Notice that
 \begin{align}
\left\|\widetilde{U}_k-U_*\right\|^2_{\rm F}&=\left\|\widetilde{U}_k-\mathbf{V}\!_k\mathbf{V}_k^TU_*-
(I-\mathbf{V}\!_k\mathbf{V}_k^T)U_*\right\|^2_{\rm F}\nonumber\\
&=\left\|\widetilde{U}_k-\mathbf{V}\!_k\mathbf{V}_k^TU_*\right\|^2_{\rm F}+\left\|(I-\mathbf{V}\!_k\mathbf{V}_k^T)U_*\right\|^2_{\rm F}\nonumber\\
&=\left\|\widetilde{U}_k-\mathbf{V}\!_k\mathbf{V}_k^TU_*\right\|^2_{\rm F}+\e_k^2.
\label{eq1802}
\end{align}
Next, we consider $\left\|\widetilde{U}_k-\mathbf{V}\!_k\mathbf{V}_k^TU_*\right\|^2_{\rm F}$. We have that
\begin{align*}
\left\|\widetilde{U}_k-\mathbf{V}\!_k\mathbf{V}_k^TU_*\right\|^2_{\rm F}
=&\|\widetilde{U}_k\|^2_{\rm F}+\|\mathbf{V}\!_k\mathbf{V}_k^TU_*\|^2_{\rm F}-2~{\rm tr}(\widetilde{U}_k^T\mathbf{V}\!_k\mathbf{V}_k^TU_*)\\
=&\|U_*\|^2_{\rm F} +\|\mathbf{V}\!_k\mathbf{V}_k^TU_*\|^2_{\rm F}-2~{\rm tr}(M)\\
=&\|U_*\|^2_{\rm F}-\|\mathbf{V}\!_k\mathbf{V}_k^TU_*\|_{\rm F}^2+2~\big(\|\mathbf{V}\!_k\mathbf{V}_k^TU_*\|^2_{\rm F}-{\rm tr}(M)\big)\\
=&\left\|(I-\mathbf{V}\!_k\mathbf{V}_k^T)U_*\right\|^2_{\rm F}+2~[{\rm tr}(M^TM)-{\rm tr}(M)]\\
=&\e_k^2+2~[{\rm tr}(M^TM)-{\rm tr}(M)].
\end{align*}
From $M=M^T\succ \bf O$ and $\sigma_i(M)=\sigma_i(\mathbf{V}\!_k\!\mathbf{V}_k^TU_*)\leq 1,~i=1,2,\ldots,\ell$, we get
$$
{\rm tr}(M^TM)-{\rm tr}(M)=\sum_{i=1}^{\ell}\sigma_i^2(M)-\sum_{i=1}^{\ell}\sigma_i(M)
=\sum_{i=1}^{\ell}\Big(\sigma_i^2(M)-\sigma_i(M)\Big)\leq 0,
$$
and $\big\|\widetilde{U}_k-\mathbf{V}\!_k\mathbf{V}_k^TU_*\big\|_{\rm F}\!\leq\e_k$.
A combination of \eqref{eq21.57} and \eqref{eq1802} yields \eqref{eq1724}.
\end{proof}
\begin{remark}
We note that
$$
\mu_1+\gamma_1\leq \|H\|_2+\|\Lambda_*\|_2\overset{\eqref{th1}}{=}\|H\|_2+\|HU_*+G\|_2\leq 2\|H\|_2+\|G\|_2.
$$
That is, $(\mu_1+\gamma_1)$ is uniformly bounded, and $\e_k\rightarrow 0$ implies $f(U_k)\rightarrow f(U_*)$.
\end{remark}

$\bullet$ Third, we show the convergence of  $U_k$. To do this, we pay special attention to the distance between the global optimal solution $U_*$ and the approximate solution $U_k$.

Notice that $U_*$ may be non-unique. In this case, the convergence of $U_k$ is difficult to
define. We first establish a sufficient condition for the uniqueness of $U_*$.
\begin{theorem}\label{thm2023}
Let $U_*$ be a global optimal solution of \eqref{1}, and define
\begin{equation}\label{eq1325}
\delta(U_*):=\inf_{X\in  U_*+\mathbb{O}^{n\times\ell}\atop X \neq\bf O}\frac{{\rm tr}(X^THX)+{\rm tr}(X\Lambda_*X^T)}{\|X\|^2_{\rm F}},
\end{equation}
where
$$
U_*+\mathbb{O}^{n\times\ell}=\big\{X~|~X=U_*+Q~~where~~Q\in \mathbb{O}^{n\times\ell}  \big\}.
$$
Then we have that
\begin{itemize}
\item[(i)]$\delta(U_*)\geq0$. Moreover, if $\delta(U_*)>0$, then $U_*$ is the unique global optimal solution to \eqref{1}.
\item[(ii)]
If the infimum in \eqref{eq1325} is attainable, then $\delta(U_*)>0$ if and only if $U_*$ is a unique global optimal solution to \eqref{1}.
\item[(iii)] We have $\delta(U_*)\geq\la_{\min}\big((I_{\ell}\otimes H)\!+\!(\Lambda_*\otimes I_n)\big)=\mu_n+\gamma_\ell$.
Specifically, if $\ell\!=\!1$, then $\Lambda_*\in\mathbb{R}$ is a scalar, and $\delta(U\!_*)=\la_{\min}(H)+\Lambda_*=\!\mu_n\!+\!\Lambda_*\geq \!0$.
\end{itemize}
\end{theorem}
\begin{proof}
(i) We prove it by contradiction. Suppose that $\delta(U_*)<0$, there is  a
matrix $\widetilde{X}\in( U_*+\mathbb{O}^{n\times\ell})\backslash \{\bf O\}$, such that
$$
\frac{{\rm tr}(\widetilde{X}^TH\widetilde{X})+{\rm tr}(\widetilde{X}\Lambda_*\widetilde{X}^T)}{\|\widetilde{X}\|^2_{\rm F}}<0.
$$
Hence, there is a matrix $L\in\mathbb{O}^{n\times\ell}$, such that $\widetilde{X}=U_*+L$. By \eqref{eq21.54},
$$0\leq \frac{f(-L)-f(U_*)}{\|U_*+L\|^2_{\rm F}}=\frac{{\rm tr}(\widetilde{X}^TH\widetilde{X})+{\rm tr}(\widetilde{X}\Lambda_*\widetilde{X}^T)}{\|\widetilde{X}\|^2_{\rm F}}<0,$$
which is a contradiction. As a result, we have $\delta(U_*)\geq0$.

Moreover, if $\delta(U_*)>0$ while $U_*$ is non-unique, then there is a matrix $U_{**}\in\mathbb{O}^{n\times\ell}$ such that $U_*\neq U_{**}$ and
 $f(U_*)=f(U_{**})$. As $ U_*-U_{**}\in( U_*+\mathbb{O}^{n\times\ell})\backslash \{\bf O\}$, it follows that
 \begin{align*}
 0<\delta(U_*)&\leq \frac{{\rm tr}[(U_*-U_{**})^TH(U_*-U_{**})]+{\rm tr}[(U_*-U_{**})\Lambda_*(U_*-U_{**})^T]}{\|U_*-U_{**}\|^2_{\rm F}}\\
              &\overset{\eqref{eq21.54}}{=} \frac{f(U_{**})-f(U_*)}{\|U_*-U_{**}\|^2_{\rm F}}=0,
 \end{align*}
 a contradiction. Thus, $U_*$ is a unique global optimal solution to \eqref{1}.

(ii) If the infimum in \eqref{eq1325} is attainable, then
 there is a matrix  $X_*\in( U_*+\mathbb{O}^{n\times\ell})\backslash \{\bf O\}$, such that
$$
X_*=\arg\min_{X\in U_*+\mathbb{O}^{n\times\ell}\atop X \neq\bm O}\frac{{\rm tr}(X^THX)+{\rm tr}(X\Lambda_*X^T)}{\|X\|^2_{\rm F}}.
$$
Thus, there is a matrix $\widehat{X}_{*}\in \mathbb{O}^{n\times\ell}$ such that $X_*=U_*-\widehat{X}_{*}$.
So we obtain from \eqref{eq21.54} that
 $$
f(\widehat{X}_{*})\!-\!f(U_*)\!=\!{\rm tr}(X_*^T\!HX_*)\!+\!{\rm tr}(X_*\Lambda_*X_*^T)\!=\!\delta(U_*)\cdot\|X_*\|^2_{\rm F}.
$$
If $U_*$ is the unique global optimal solution to \eqref{1}, then
$
f(\widehat{X}_{*})-f(U_*)>0
$.
That is, $\delta(U_*)>0$, and we have (ii) from (i).

(iii) We have that
\begin{align}
\delta(U_*)
\!\geq&\!\inf_{X\in \mathbb{R}^{n\times \ell}\atop X \neq\bm O}\!\frac{{\rm tr}(X^THX)\!+\!{\rm tr}(X\Lambda_*X^T)}{\|X\|^2_{\rm F}}
=\inf_{X\in \mathbb{R}^{n\times \ell}\atop X \neq\bm O}\!\frac{{\rm tr}(X^THX)\!+\!{\rm tr}(\Lambda_*X^T\!X)}{\|X\|^2_{\rm F}}
\nonumber\\
\overset{\eqref{eq1148}}{=}&\!\inf_{X\in \mathbb{R}^{n\times \ell}\atop X \neq\bm O}\frac{({\rm vec}(X))^T\cdot\Big((I_{\ell}\otimes H)\!+\!(\Lambda_*\otimes I_n)\Big)\cdot{\rm vec}(X)}{\|{\rm vec}(X)\|^2_{\rm 2}}\nonumber\\
\overset{\eqref{eq2031}}{=}&\la_{\min}(\mathcal{H}_*)
\overset{\eqref{eq2031b}}{=}\mu_n+\gamma_{\ell}.\nonumber
\end{align}

In particular, if $\ell\!=\!1$, then \eqref{1} reduces to a trust-region subproblem. In this case, $\Lambda_*\!=\!\la_*\!\in\!\mathbb{R}$ is a scalar, $H_*:=H\!+\!\la_* I\!\succcurlyeq\!\bf O$ \cite[Lemma 2.1]{W.W},
$\mathbb{O}^{n\times \ell}\!\!=\mathcal{B}:=\{\bm x~\big|~\|\bm x\|_2\!=\!1\}$, and
$$
\delta(U_*)\!=\!\inf_{\bm x\in \mathcal{S}} \bm x^TH_*\bm x,
~~~
{\rm with}
~~~
\mathcal{S}=\left\{\bm x~\Big|~\bm x=\frac{U_*-\bm p}{\|U_*-\bm p\|},~~{\rm where}~~U_*\neq\bm p\in\mathcal{B}\right\}\subseteq\mathcal{B}.
$$

Denote by $\widetilde{\mathcal{S}}=\left\{\bm x~\Big|~\bm x=\frac{\bm p -U_*}{\|\bm p -U_*\|},~~{\rm where}~~U_*\neq\bm p\in\mathcal{B}\right\}$, we see
that
$$
\inf_{\bm x\in\mathcal{S}\cup\widetilde{\mathcal{S}}}\bm x^TH_*\bm x
=\min\left\{\inf_{\bm x\in\mathcal{S}}\bm x^TH_*\bm x,~\inf_{\bm x\in\widetilde{\mathcal{S}}}\bm x^TH_*\bm x\right\},
$$
and
$
\inf_{\bm x\in\mathcal{S}}\bm x^TH_*\bm x=\inf_{\bm x\in\widetilde{\mathcal{S}}}\bm x^TH_*\bm x.
$
Therefore, $\delta(U_*)=\underset{\bm x\in \mathcal{S}\cup\widetilde{\mathcal{S}}}{\inf} \bm x^TH_*\bm x$.
As
$$
\mathcal{S}\cup\widetilde{\mathcal{S}}=\left\{\bm x~\Big|~\bm x=\pm\frac{\bm p-U_*}{\|\bm p-U_*\|},~~{\rm where}~~U_*\neq\bm p\in\mathcal{B}\right\}
$$
is dense on $\mathcal{B}$ \cite[p. 91]{7}, it follows that
$$
\delta(U_*)=\inf_{\bm x\in \mathcal{B}} \bm x^T(H\!+\!\la_*I)\bm x
=\min_{\bm x\in \mathcal{B}} \bm x^T(H\!+\!\la_*I)\bm x=\mu_n+\la_*\geq 0,
$$
which completes the proof.
\end{proof}

We are ready to consider the convergence of $U_k$.
\begin{theorem}\label{Thm3.7}
Let $U_*$ be the global optimal solution of \eqref{1}.
If $\delta(U_*)>0$ and $\|(I-\mathbf{V}\!_k\mathbf{V}_k^T)U_*\|_2<1$, then we have
\begin{equation}\label{eq1725}
\|U_k-U_*\|_{\rm F}\leq \sqrt{\frac{2(\mu_1+\gamma_1)}
{\delta(U_*)}}\cdot\e_k.
\end{equation}
Specifically, if $\mathcal{H}_*\succ\bf O$, then
\begin{equation}\label{eq17.25}
\|U_k-U_*\|_{\rm F}\leq \sqrt{2\varkappa_*}\cdot\e_k,
\end{equation}
where $\varkappa_*$ is the 2-condition number of $\mathcal{H}_*$.
\end{theorem}
\begin{proof}
We notice that $U_*-U_{k}\in ( U_*+\mathbb{O}^{n\times\ell})\cup(-U_*+\mathbb{O}^{n\times\ell})$.
It follows from \eqref{eq21.54} and \eqref{eq1325} that
\begin{align}
f(U_k)\!-\!f(U_*)
             &\!=\!\frac{{\rm tr}[(U_*\!-\!U_k)^T\!H(U_*-U_k)]+{\rm tr}[(U_*\!-\!U_k)\Lambda_* (U_*\!-\!U_k)^T]}{\|U_*\!-\!U_k\|^2_{\rm F}}\cdot\|U_*\!-\!U_k\|^2_{\rm F}\nonumber\\
             &\geq \delta(U_*)\cdot\|U_*-U_k\|^2_{\rm F},\label{eq21.28}
\end{align}
so we obtain \eqref{eq1725} from combining \eqref{eq21.28} and \eqref{eq1724}.

If $\mathcal{H}_*\succ\bf O$, it is seen from Theorem \ref{thm2023} (iii) and \eqref{eq21.28} that
\begin{align*}
\small
\|U_*\!-\!U_k\|^2_{\rm F}\leq&\frac{\delta(U_*)\|U_*\!-\!U_k\|^2_{\rm F}}{\mu_n\!+\!\gamma_{\ell}}\leq
\frac{f(U_k)\!-\!f(U_*)}{\mu_n\!+\!\gamma_{\ell}}\\
\overset{\eqref{eq1724}}{\leq}& \frac{2(\mu_1\!+\!\gamma_1)}{\mu_n\!+\!\gamma_{\ell}}\cdot\e^2_k
\!\overset{\eqref{eq2031}}{=}\!2\varkappa\cdot\e^2_k,
\end{align*}
which yields \eqref{eq17.25}.
\end{proof}
\begin{remark}
Theorem \ref{Thm3.7} indicates that $\delta(U_*)$ plays an important role in the convergence of $U_k$. More precisely, $U_k$ may converge slowly as $\delta(U_*)$ is close to zero. Specifically, $U_k$ is difficult to define the convergence as $\delta(U_*)=0$, which coincides with the results established in Theorem \ref{thm2023}.
\end{remark}

$\bullet$ Fourth, we consider the KKT error $\|HU_k+U_k\Lambda_k+G\|_{\rm F}$
and the upper bound on $\|\Lambda_*-\Lambda_k\|_{\rm F}$, where $\Lambda_k$ is defined in \eqref{eq10}.
\begin{theorem}
Denote by
$$
R_k= HU_k+U_k\Lambda_k+G.
$$
Then
\begin{itemize}
\item[(i)]
We have that
 \begin{equation}\label{eq22}
\max\big\{\|R_k\|_{\rm F},\|\Lambda_*-\Lambda_k\|_{\rm F}\big\}\leq \|\mathcal{H}_*\|_2\cdot\|U_k-U_*\|_{\rm F}.
\end{equation}
\item[(ii)]
If $\mathcal{H}_*\succ\bf O$ and $\|(I-\mathbf{V}\!_k\mathbf{V}_k^T)U_*\|_2<1$, then
\begin{equation}\label{eq2001}
\max\big\{\|R_k\|_{\rm F},\|\Lambda_*-\Lambda_k\|_{\rm F}\big\}\leq
\sqrt{2}\|\mathcal{H}_*\|_2\cdot\e_k.
\end{equation}
\item[(iii)]
If  $\mathcal{K}_q$ is an invariant subspace of $H$,  then
$$
\|R_q\|_{\rm F} = \|HU_q+U_q\Lambda_q+G\|_{\rm F}= 0~~{  and}~~ -U_q^T G\succcurlyeq\bf O.
$$
That is, $U_q$
satisfies the first order optimality in Theoem \ref{ThmKKT} and also the necessary condition \eqref{glopt} for a global minimizer.
\end{itemize}
\end{theorem}
\begin{proof}
(i) We notice that
$$
P_k =\arg\min_{P\in \mathbb{O}^{k \ell\times \ell}}\Big\{{\rm tr}(P^TT_k P)+2~{\rm tr}(P^T{G}_k)\Big\},
$$
and $G = \mathbf{V}\!_k\mathbf{V}_k^TG=\mathbf{V}\!_kG_k$. From \eqref{eq1653},
 \begin{align}
R_k=& HU_k+U_k \Lambda_k+G\nonumber\\
  =&H\mathbf{V}\!_kP_k+\mathbf{V}\!_kP_k\Lambda_k+G\nonumber\\
  =&\left(\mathbf{V}\!_kT_k+W_{k+1}N_{k}(E^{(k\ell)}_{\ell})^T\right)
  P_k+\mathbf{V}\!_kP_k\Lambda_k+\mathbf{V}\!_kG_k\nonumber\\
  =&\mathbf{V}\!_k\big(T_kP_k+P_k\Lambda_k+G_k\big)
  +W_{k+1}N_{k}(E^{(k\ell)}_{r})^TP_k\label{eq1633}\\
  \overset{\eqref{eq10}}{=}&W_{k+1}N_{k}(E^{(k\ell)}_{\ell})^TP_k,\nonumber
 \end{align}
and $\mathbf{V}_k^T R_k=\bf O$. Thus, we have from $U_k=\mathbf{V}\!_k P_k$ that
\begin{align*}
\|HU_k+U_k\Lambda_*+G\|^2_{\rm F} =\| R_k+U_k(\Lambda_*-\Lambda_k)\|^2_{\rm F}
                                 =\|R_k\|^2_{\rm F}+\|\Lambda_*-\Lambda_k\|^2_{\rm F}.
\end{align*}
As a result,
\begin{equation}\label{eq1600}
\max\big\{\|R_k\|_{\rm F}, \|\Lambda_*-\Lambda_k\|_{\rm F}\big\}\leq
\|HU_k+U_k\Lambda_*+G\|_{\rm F}.
\end{equation}
Notice that
\begin{align}
\|HU_k+U_k\Lambda_*+G\|_{\rm F}\overset{\eqref{th1}}{=}&\|HU_k+U_k\Lambda_*-(HU_*+U_*\Lambda_*)\|_{\rm F}\nonumber\\
=&\|H(U_k-U_*)+(U_k-U_*)\Lambda_*\|_{\rm F}\nonumber\\
=&{\Big\|}{\rm vec}\big(H(U_k-U_*)+(U_k-U_*)\Lambda_*\big) {\Big\|}_{2}
\nonumber\\
=&\|\left[(I_{\ell}\otimes H)+(\Lambda_*\otimes I_{n})\right]{\rm vec}(U_k-U_*)\|_2\label{eq15.42}\\
\leq&\|\mathcal{H}_*\|_2\cdot\|U_k-U_*\|_{\rm F}.\nonumber
\end{align}
So we have \eqref{eq22} from \eqref{eq1600} and
\eqref{eq15.42}.

(ii)
If $\mathcal{H}_*=(I_{\ell}\otimes H)+(\Lambda_*\otimes I_n)\succ\bf O$, then
\begin{align*}
&\big\|\big[(I\otimes H)+(\Lambda_*\otimes I)\big]{\rm vec}(U_k-U_*)\big\|^2_2\\
\leq&\big\|\mathcal{H}_*^{\frac{1}{2}}\big\|^2_2\cdot\big\|\mathcal{H}_*^{\frac{1}{2}}{\rm vec}(U_k-U_*)\big\|^2_2\\
=&{\|\mathcal{H}_*\|_2}\cdot
{\big[{\rm vec}(U_k-U_*) \big]^T (I\otimes H+\Lambda_*\otimes I)~{\rm vec}(U_k-U_*)}\\
=&{\|\mathcal{H}_*\|_2}\cdot\left(
{{\rm tr}[(U_k\!-\!U_*)^TH(U_k\!-\!U_*)]\!+\!{\rm tr}[\Lambda_* (U_k\!-\!U_*)^T(U_k\!-\!U_*)
]} \right)\\
\overset{\eqref{eq21.54}}{=} & \|\mathcal{H}_*\|_2\cdot  \big(f(U_k)-f(U_*)\big)
\overset{\eqref{eq1724}}{\leq}2\|\mathcal{H}_*\|_2(\mu_1+\gamma_1)\cdot\e_k^2
\\
\leq&2\|\mathcal{H}_*\|^2_2\cdot\e_k^2 .
\end{align*}
A combination of \eqref{eq1600} and \eqref{eq15.42} gives \eqref{eq2001}.

(iii) In this case, we have from \eqref{eq1653} that $H\mathbf{V}\!_q=\mathbf{V}\!_qT_q$.
Recall that $G=\mathbf{V}\!_q\mathbf{V}_q^TG=\mathbf{V}\!_qG_q$. Hence,
\begin{align*}
R_q = HU_q+U_q\Lambda_q+G=H\mathbf{V}\!_qP_q+\mathbf{V}\!_q(P_q\Lambda_q+G_q)
\overset{\eqref{eq10}}{=}(H\mathbf{V}\!_q-\mathbf{V}\!_q T_q)P_q=\bf O.
\end{align*}
By \eqref{glopt}, $-P_q^TG_q\succcurlyeq\bf O$, and $-U_q^TG=-P_q^T\mathbf{V}_q^T\mathbf{V}_qG_q=-P_q^TG_q\succcurlyeq\bf O$.
\end{proof}

\begin{remark}
It is known that if the block Lanczos process terminates at the $q$-th step, then $\mathcal{K}_q$ is an invariant subspace of $H$, and $q$ is no larger than the number of distinct eigenvalues of $H$ \cite{2,Y.S,Stewart}. This can happen in some applications such as the orthogonal least squares regression (OLSR) model \cite{OLSR1} for supervised learning, where $H$ is often a low-rank matrix.
\end{remark}

\section{Numerical experiments}\label{Sec4}
In this section, we perform numerical experiments to illustrate the numerical behavior of the proposed method.
All the numerical experiments were run on a AMD R7 5800H CPU 3.20 GHz with 16GB RAM under Windows 11 operation system.
The experimental results are obtained from using MATLAB R2022a implementation with machine precision $u_{\rm machine} \approx 2.22 \times 10^{-16}$.
To show the efficiency of Algorithm \ref{alg}, we compare it with seven state-of-the-art approaches for solving \eqref{1}, including
the RTR method \cite{Abs}, the SCFRTR method \cite{Z.L}, the PCAL method \cite{Yuan4},
 the GP-BB method \cite{Yuan1}, the WYBB method \cite{Wen}, the JDCP method \cite{Jiang}, as well as the GPI method \cite{Niee}.

In all the experiments, we first normalize $H$ and $G$ by using $\|G\|_{\rm F}$,
and use the following stopping criterion \cite{Jiang,Wen,Z.L}
\begin{equation}\label{eq1234}
\frac{|f(U_k)-f(U_{k+1})|}{|f(U_k)|\!+\!1}\leq \varepsilon_f,~\frac{\|U_k-U_{k+1}\|_{\rm F}}{\sqrt{n}}\leq \varepsilon_U,~\|R_k\|_{\rm F}\leq \varepsilon_{g},~{\rm and}~k\leq k_{\max},
\end{equation}
with $\varepsilon_f= 10^{-10}$, $\varepsilon_U=10^{-6}$, $\varepsilon_{g}= 10^{-5}$ and $k_{\max}=1000$.

\begin{figure}[H]\caption{Example 4.1: Numerical results on the synthetic data.}\label{fig:1114}
\centering
{\includegraphics[width=7.2cm,height=1.25cm]{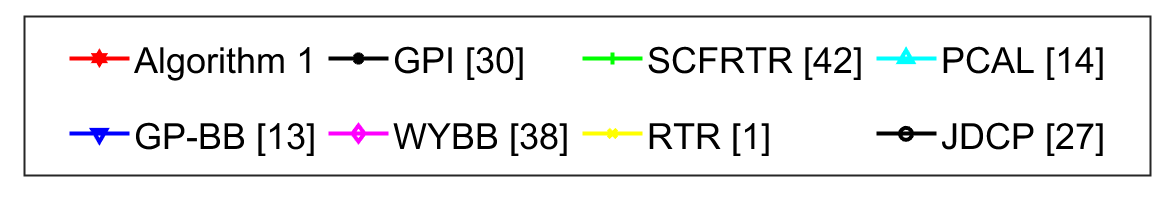}}\\
\subfigure[$\ell=10$: The relative objective function difference $f_{\rm err}^{\rm (rel)}$.]{\includegraphics[width=5.123cm,height=3.76cm]{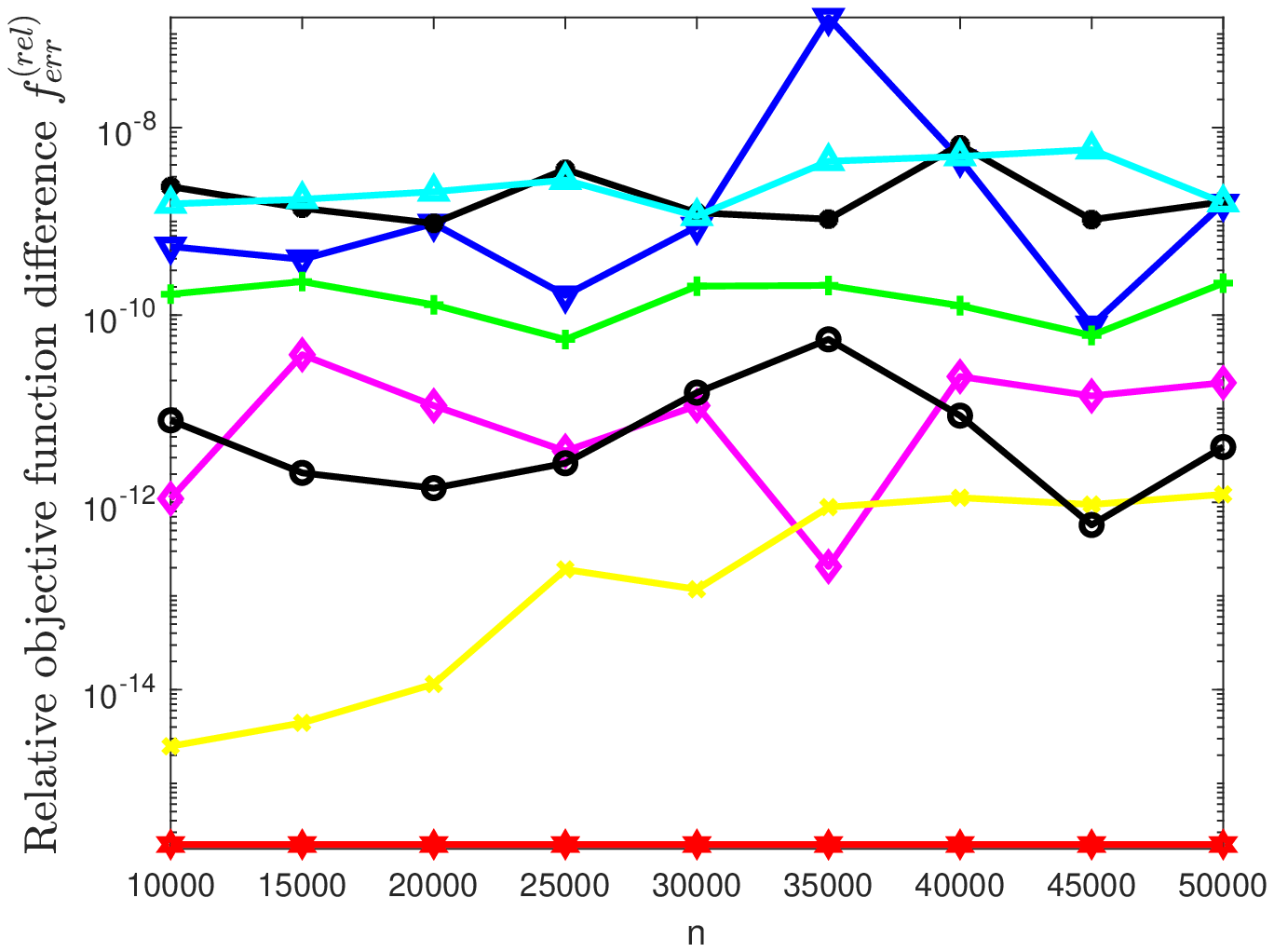}}
\subfigure[$\ell=10$: KKT error]{\includegraphics[width=5.123cm,height=3.76cm]{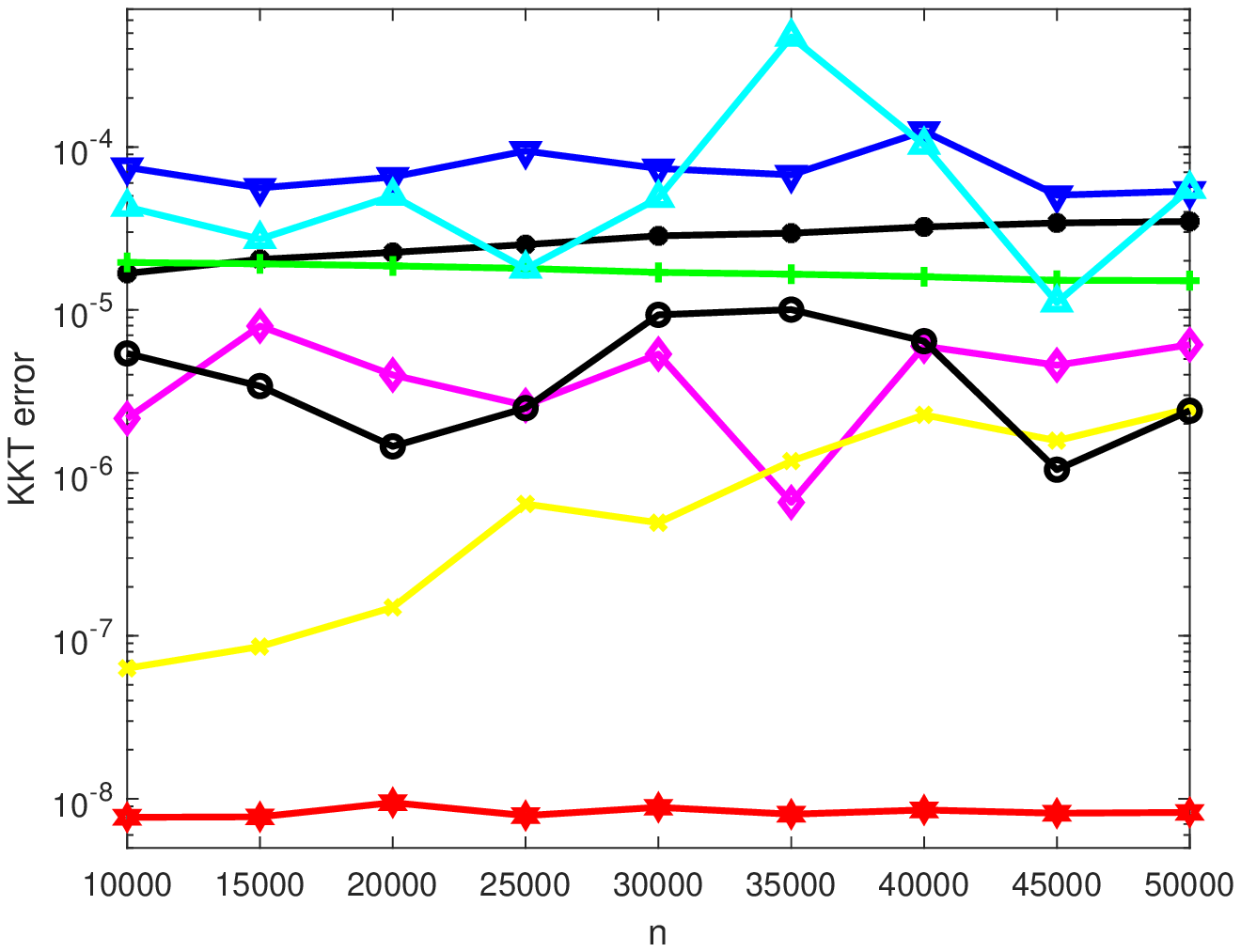}}
\subfigure[$\ell=10$: CPU time]{\includegraphics[width=5.123cm,height=3.76cm]{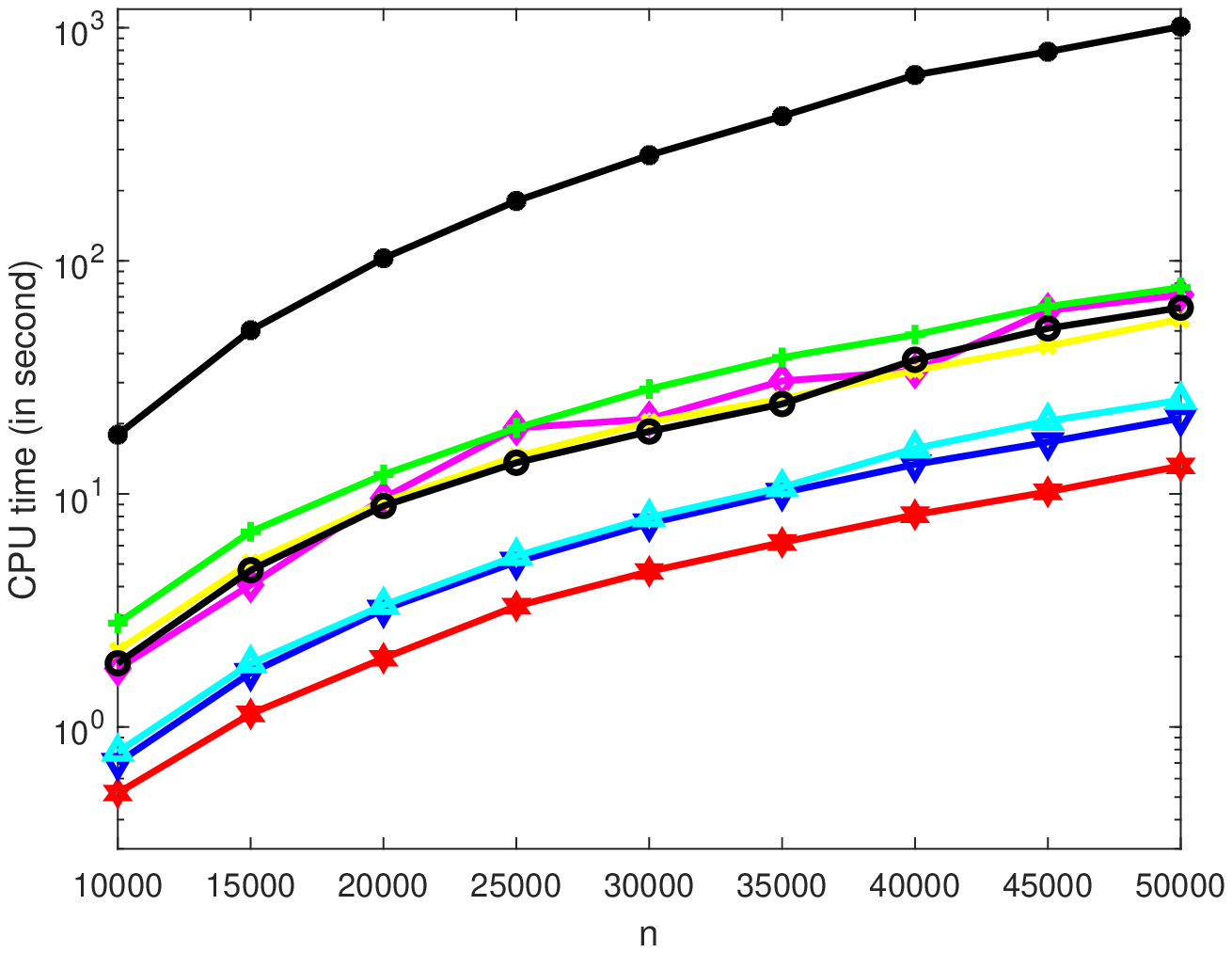}}
\subfigure[$\ell=20$: The relative objective function difference $f_{\rm err}^{\rm (rel)}$.]{\includegraphics[width=5.123cm,height=3.76cm]{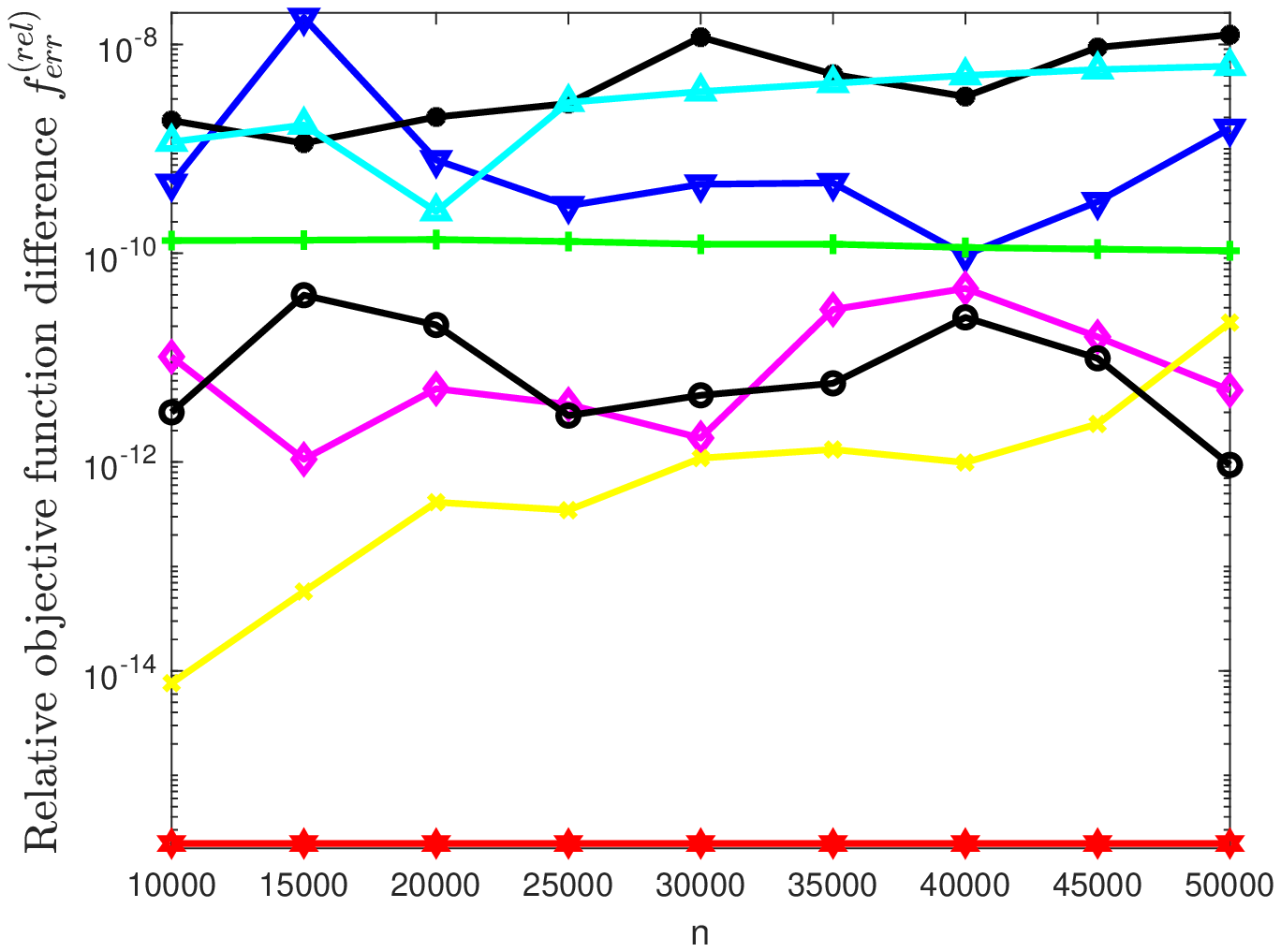}}
\subfigure[$\ell=20$: KKT error]{\includegraphics[width=5.123cm,height=3.76cm]{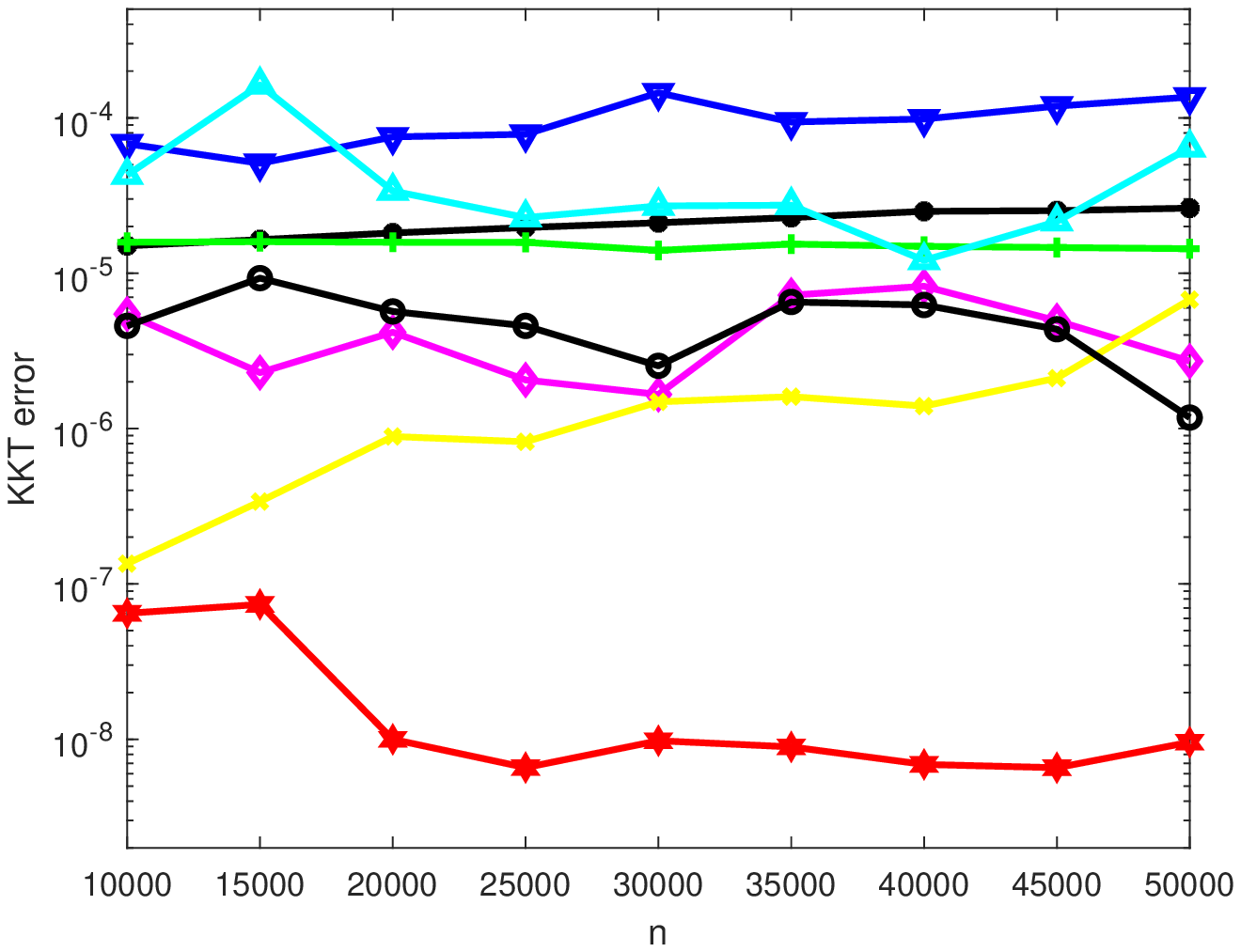}}
\subfigure[$\ell=20$: CPU time]{\includegraphics[width=5.123cm,height=3.76cm]{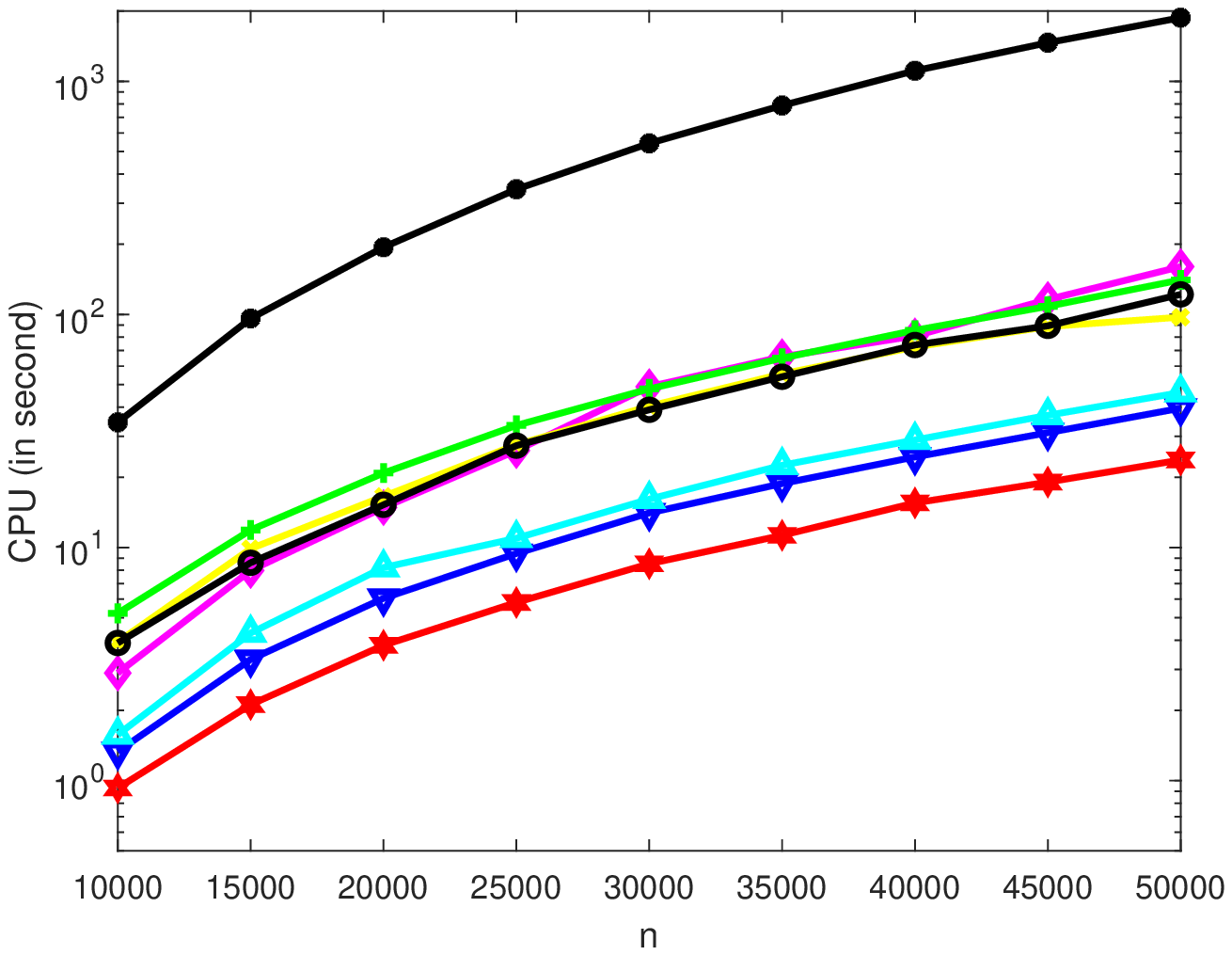}}
\end{figure}

In the block Lanczos process, we make use of full reorthogonalization process when necessary.
We stress that
an advantage of the proposed method is that one can compute the KKT error cheaply. More precisely, we have from \eqref{eq1633} that
\begin{align}
\|R_k\|_{\rm F}&=\|HU_k+U_k \Lambda_k+G\|_{\rm F}\nonumber\\
               &=\sqrt{ \left\|T_kP_k+P_k\Lambda_k+G_k\right\|^2_{\rm F}+\left\|N_{k}(E^{(k\ell)}_{\ell})^TP_k\right\|^2_{\rm F}}\nonumber\\
               &=\sqrt{ \left\|T_kP_k\!+\!P_k\Lambda_k\!+\!G_k\right\|^2_{\rm F}+\left\|N_{k}P_k\big((k\!-\!1)\ell\!+\!1:k\ell,:\big)\right\|^2_{\rm F}} .
\end{align}
Moreover,  $f(U_k)=\widetilde{f}(P_k)$ and $f(U_{k+1})=\widetilde{f}(P_{k+1})$; refer to \eqref{eqn2.5}. As $k$ increases, the main overhead in each iteration of our method lies in solving \eqref{3} by using the RTR method. Thus, we solve \eqref{3} every 5 steps in practical calculations.

To measure the accuracy of the approximations from the algorithms, we make use of the relative objective function difference defined as \cite{3}
\begin{equation}
f^{\rm (rel)}_{\rm err} := \frac{f(\breve{{U}}_*)-f(U_{\rm best})}{|f(U_{\rm best})|},
\end{equation}
where ${\breve{U}}_*$ is the computed solution of each method and $U_{\rm best}$ denotes the solution with the smallest objective value among all the solvers. Thus, $f^{\rm (rel)}_{\rm err}=0$ means that $\breve{{U}}_*=U_{\rm best}$.

\subsection{Test on synthetic data}\label{sec4.1}
In this subsection,  we make experiments on some synthetic data generated by using the MATLAB built-in function {\tt sprand}:
$$
H = B+B^T,~~{\rm where}~~B={\tt sprand(n,n,density)},~~G={\tt randn(n,\ell)},
$$
where $density\!=\!0.05$, $n\!=\!10000,15000,20000,\ldots,50000$, and $\ell\!=\!10,20$, respectively.

The numerical results of the eight algorithms are ploted in Figure \ref{fig:1114}. It is seen from the figure that
both the relative objective function difference $f_{\rm err}^{\rm (rel)}$ and the KKT errors of of Algorithm \ref{alg} are the smallest, and our algorithm is the fastest one among the eight algorithms.





\subsection{Test on the orthogonal least squares regression for feature extraction}\label{sec2}
Orthogonal least squares regression (OLSR) is a popular supervised learning method for linear discriminant analysis (LDA) \cite{OLSR1}.  Let $A = [\bm a_1,\ldots,\bm a_m]\in \mathbb{R}^{n\times m}$ be the whole database with $\ell$ classes, where $m$ is the number of samples and $n$ is the number of
features.
Let $\widehat{A}=[\hat{\bm a}_1,\ldots,\hat{\bm a}_{\tilde{m}}]\in\mathbb{R}^{n\times\tilde{m}}$ be training data set, and
$B = [\bm b_1,\ldots,\bm b_{\tilde{m}}] \in \mathbb{R}^{\ell \times \tilde{m}}$ be the corresponding
class indicator matrix, where $\tilde{m}$ is the number of training samples,
and $\bm b_i = \bm e_j\in \mathbb{R}^{\ell}$ if  the sample $\hat{\bm a}_i$ is in the $j$-th class, $1\leq i\leq\tilde{m},1\leq j\leq\ell$, where $\bm e_j$ is the $j$-th column of the identity matrix. In the experiment, we randomly choose 30\% of the total samples as the training set. The details of of the fourteen data sets are listed in Table \ref{11}.

Let $\widetilde{A},\widetilde{B}$ be the centered matrix of $\widehat{A},\widehat{B}$, respectively.
In the OLSR method, one aims to seek $U \in \mathbb{O}^{ n\times \ell}$  such that \cite{OLSR1}


\begin{equation}\label{OLSR}
\min_{U \in \mathbb{O}^{ n\!\times\!\ell}} \left\{{\rm tr}(U^T\!HU)\!+\!2{\rm tr}(U^T\!G)\right\},
\end{equation}
where $H\!=\!\widetilde{A}^T\widetilde{A}\!\in\!\mathbb{R}^{n\times n}$
and $G \!=\!\widetilde{A}^T\widetilde{B}\!\in\!\mathbb{R}^{n\times\ell}$.
We run the eight algorithms on the fourteen databases, and the numerical results are reported in Table \ref{tab1} and Table \ref{tab2}.
Specifically, if the CPU time of an algorithm exceeds 3600 seconds or the KKT error $\|H\breve{U}_*+\breve{U}_* \breve{\Lambda}_*+G\|_{\rm F}\geq 1$,
we declare that the algorithm fails to converge  and denote it by ``--"".

We observe from Table \ref{tab1} and \ref{tab2} that Algorithm \ref{alg} is more powerful than the other seven popular algorithms for solving the OLSR model \eqref{OLSR}.
More precisely, Algorithm \ref{alg} is the best in terms of the values of $f^{\rm (rel)}_{\rm err}$, and the KKT errors from Algorithm \ref{alg} is the smallest except for the {\tt ORL} and {\tt Text-1} databases. Indeed, the KKT errors of our algorithm is about two to five orders lower than those of the others.
Moreover, our algorithm is the fastest one except for the {\tt YouTubeFace} database, which ours is the second fastest one.

 \begin{table}[H]
\begin{center}\caption{Summary of test datasets in Example \ref{sec2}.}\label{11}
\begin{tabular}{l||c|c|c|c}
\hline \hline
 Datasets &  \!Feature ($n$)\! & \tabincell{c}{\!Number of\\ \!samples ($m$)}&
  \tabincell{c}{\!Number of\\ \!classes ($\ell$)} & Background \\
\hline
\tt ORL\tablefootnote{\url{http://featureselection.asu.edu/datasets.php}.} &10304 & 400 &40 &\multirow{3}{*}{Image} \\
\tt Yale\tablefootnote{\url{https://www.face-rec.org/databases/}.}  & 10000 & 165 & 15&\\
\tt YouTubeFace\tablefootnote{\url{https://www.cs.tau.ac.il/~wolf/ytfaces/.}} &16384 & 56653   &17  \\
\hline
{\tt CLL\_SUB\_111}\tablefootnote{The databases
{\tt CLL\_SUB\_111, SMK\_CAN\_187,  GLI\_85, leukemia} and
{\tt nci9}  are available at
https://jundongl.github.io/scikit-feature/datasets.html.} &11340	&111&	3&\multirow{5}{*}{Biological}\\

\tt SMK\_CAN\_187 &19993 &187	&2\\

\tt GLI\_85 &22283	&85&	2\\

\tt leukemia  & 7070& 72&   2  \\

\tt nci9  &9712 &60&	9\\
\hline
\tt 20Newsgroups\tablefootnote{The databases {\tt 20Newsgroups} and {\tt RCV1\_4Class.} are available at \url{http://www.cad.zju.edu.cn/home/dengcai/Data/TextData.html}} &  26214   &  18846   &20 & \multirow{7}{*}{Text}\\
\tt RCV1\_4Class &29992  & 9625 &4 \\
\tt Text-1\tablefootnote{The databases
 {\tt Text-1, Cora-HA, Cora-OS} and {\tt Core-PL}  are available at \url{http://www.escience.cn/people/fpnie/papers.html}.}      &7511  &1946 &2\\
\tt Cora-HA  & 3989 & 400 & 7\\
\tt Cora-OS&6737 & 1246 & 4\\
\tt Core-PL & 7949 & 1575 & 9 \\
\hline\hline
\end{tabular}
\end{center}
\end{table}
\begin{table}[H]
\scriptsize
  \centering
  \caption{Example \ref{sec2}: Numerical experiments on the OLSR model \eqref{OLSR} for some data sets in Table \ref{11}, where the best results are in bold.}\label{tab1}
    \begin{tabular}{p{2cm}|p{2.4cm}|p{1.6cm}|p{1.3cm}}
     \hline
      Datasets   &CPU(s)&KKT error&~~$f^{(\rm rel)}_{\rm err}$ \\
    \hline
    \multirow{8}{*}{\tt ORL}&
    {\bf Alg.\ref{alg}}: {\bf 1.01} \newline{}
    SCFRTR: 10.62\newline{}
    WYBB: 31.11\newline{}
    JDCP: 32.05\newline{}
    RTR: 24.25\newline{}
    GP-BB: 6.85\newline{}
    PCAL: 8.34\newline{}
    GPI: 1.22e+02
    & 8.08e-06\newline{}
    8.86e-06\newline{}
    4.58e-05\newline{}
    2.55e-05\newline{}
    {\bf 1.49e-07}\newline{}
    1.74e-03 \newline{}
    1.32e-03 \newline{}
    5.30e-03
    & 4.37e-13\newline{}1.07e-11\newline{}6.31e-08\newline{}
    8.21e-09\newline{}{\bf 0}       \newline{}1.19e-05\newline{}  4.76e-05\newline{}9.70e-05 \\ \hline

    \multirow{8}{*}{\tt Yale}
    & {\bf Alg.\ref{alg}: 1.34}\newline{}SCFRTR: 18.16\newline{}
    WYBB: 50.29\newline{}JDCP: 45.30\newline{}
    RTR: 75.44\newline{}GP-BB: 14.83\newline{}
    PCAL: 14.00\newline{}GPI: 1.80e+02
  & {\bf 8.46e-08}\newline{}
  8.72e-06\newline{}4.93e-05\newline{}5.07e-05\newline{} 9.96e-06\newline{}1.74e-03 \newline{}8.73e-04\newline{}4.60e-03
  & {\bf 0} \newline{}  3.44e-09\newline{}1.92e-07\newline{}1.83e-07
  \newline{} 6.66e-10\newline{}8.03e-05
  \newline{} 6.20e-05\newline{}6.76e-05\\ \hline

    \multirow{8}{*}{\tt YouTubeFace}
    & {\bf Alg.\ref{alg}}: 1.60e+02\newline{}SCFRTR: 8.15e+02\newline{}
    WYBB: 7.81e+02 \newline{}JDCP: 1.00e+03 \newline{}
    RTR: --\newline{}GP-BB: 41.07\newline{}
    PCAL: {\bf 29.08}\newline{}GPI: 3.31e+02
  & {\bf 8.46e-06}\newline{}
  8.72e-06\newline{}1.10e-03\newline{}1.70e-03\newline{}--
  \newline{}4.85e-03 \newline{}6.43e-03 \newline{}9.76e-03
  &{\bf 0}\newline{}1.98e-06\newline{}1.24e-02\newline{}5.42e-03
  \newline{}-- \newline{}5.81e-02
  \newline{} 7.80e-02\newline{}1.94e-01\\ \hline

    \multirow{8}{*}{\tt CLL\_SUB\_111} &
    {\bf Alg.\ref{alg}: 2.73}\newline{}SCFRTR: 5.34\newline{}
    WYBB: 30.42\newline{}JDCP: 26.06\newline{}
    RTR: 25.58\newline{}GP-BB: --\newline{}
    PCAL: --\newline{}GPI: --
    & {\bf 5.19e-09}\newline{}6.76e-08\newline{} 1.54e-03\newline{}2.80e-01\newline{} 5.96e-07\newline{}-- \newline{}--\newline{}--
    &{\bf 0} \newline{}  1.70e-07\newline{}1.20e-03\newline{}26.56
  \newline{} 1.83e-07\newline{}--
  \newline{} --\newline{}--\\ \hline

   \multirow{8}{*}{\tt SMK\_CAN\_187} &
    {\bf Alg.\ref{alg}}: 8.87\newline{}SCFRTR: 22.13\newline{}
    WYBB: 31.81\newline{}JDCP: 33.76\newline{}
    RTR: 25.58\newline{}GP-BB: 6.69\newline{}
    PCAL: {\bf 6.65}\newline{}GPI: 2.96e+02
    & {\bf 3.68e-09}\newline{}1.90e-07\newline{}1.26e-05\newline{}
    2.05e-05\newline{}2.79e-08\newline{}1.52e-02\newline{} 1.98e-02\newline{} 2.71e-02
    &{\bf 0} \newline{}  5.24e-13
\newline{}3.33e-09\newline{}5.50e-09
  \newline{} {\bf 0}\newline{}3.80e-03
  \newline{} 6.10e-03\newline{}8.25e-05\\ \hline

   \multirow{8}{*}{\tt  GLI\_85} &
   {\bf Alg.\ref{alg}: 5.69}\newline{}SCFRTR: 15.84\newline{}
   WYBB: 27.99\newline{}JDCP: 22.95\newline{}
   RTR: 45.51\newline{}
   GP-BB: --\newline{}PCAL: --\newline{}
   GPI: --
   &{\bf 1.80e-08}\newline{}1.00e-07\newline{}1.00e-02\newline{}
   9.00e-02\newline{}4.72e-08\newline{}--\newline{}-- \newline{}--
   & {\bf 0 }\newline{}  4.78e-07\newline{}1.30e-03\newline{}1.20e-01
  \newline{}2.17e-07\newline{}--
  \newline{}--\newline{}--\\ \hline

   \multirow{8}{*}{\tt leukemia}
   & {\bf Alg.\ref{alg}}: {\bf 0.55}\newline{}SCFRTR: 0.79\newline{}WYBB: 0.84\newline{}JDCP: 0.92\newline{}RTR: 1.51\newline{}
   GP-BB: 0.58\newline{}PCAL: 0.90\newline{}
   GPI: 6.25   & {\bf 1.48e-08}\newline{}1.34e-06\newline{}2.47e-04\newline{}1.28e-05\newline{}
   6.89e-06\newline{}5.39e-04 \newline{}4.32e-03 \newline{}5.23e-04
   &{\bf 0 }\newline{} 3.39e-12\newline{}3.42e-08\newline{}3.37e-10
  \newline{}1.10e-10\newline{}3.25e-07
  \newline{}3.55e-05\newline{}5.50e-07\\ \hline

    \multirow{8}{*}{\tt nci9}  &
    {\bf Alg.\ref{alg}: 0.51}\newline{}SCFRTR: 2.54\newline{}WYBB: 3.70\newline{}JDCP: 3.69\newline{}RTR: 5.44\newline{}GP-BB: 1.71\newline{}PCAL: 3.35\newline{}GPI: 15.04
    &{\bf 1.79e-08}\newline{}5.89e-08\newline{}6.80e-05\newline{}
    2.70e-05\newline{}1.28e-07\newline{}4.47e-03 \newline{}2.44e-03  \newline{}5.92e-04
    &{\bf 0} \newline{} {\bf 0}\newline{}1.17e-08\newline{}1.50e-09
  \newline{}{\bf 0}\newline{}2.42e-06
  \newline{}1.19e-05\newline{}1.59e-07 \\ \hline

   \multirow{8}{*}{\tt Text-1}&
   {\bf Alg.\ref{alg}: 1.20}\newline{}SCFRTR: 3.01\newline{}WYBB: 3.79\newline{}JDCP: 3.49\newline{}RTR: 28.59
    \newline{}GP-BB: 1.50\newline{}PCAL: 2.20\newline{}GPI: 41.36
    &4.72e-06\newline{}1.07e-06\newline{} 9.51e-06\newline{}9.29e-06\newline{}{\bf 7.44e-08}\newline{}
    2.87e-05\newline{}8.90e-05 \newline{}3.73e-05
    &5.92e-10\newline{} 6.60e-13\newline{}3.54e-08\newline{}1.66e-08
  \newline{}{\bf 0}\newline{}3.26e-07
  \newline{}8.14e-07\newline{}6.50e-07\\ \hline
 \end{tabular}
\end{table}%

\begin{table}[H]
\scriptsize
  \centering
  \caption{Example \ref{sec2}: Numerical experiments on the OLSR model \eqref{OLSR} for some data sets in Table \ref{11}, where the best results are in bold.}\label{tab2}%
    \begin{tabular}{p{2cm}|p{2.4cm}|p{1.3cm}|p{1.3cm}}
     \hline
      Datasets   &CPU(s)&KKT error&~~$f^{(\rm rel)}_{\rm err}$
    \\ \hline
    \multirow{8}{*}{\tt 20Newsgroups}
    & {\bf Alg.\ref{alg}: 56.03}\newline{}SCFRTR: 6.25e+02\newline{}WYBB: 1.87e+03\newline{}JDCP: 2.36e+03\newline{}RTR: --\newline{}GPBB: 1.16e+02\newline{}
    PCAL: 1.00e+02\newline{}GPI: 4.40e+02   &{\bf 8.93e-06}\newline{}2.03e-05\newline{}3.37e-03\newline{}2.70e-04\newline{}--
    \newline{}1.47e-02  \newline{}1.31e-02\newline{}5.49e-02
    &{\bf 0}\newline{}6.93e-08\newline{}3.26e-05\newline{}5.41e-06
  \newline{}--\newline{}3.29e-02
  \newline{}2.65e-02\newline{}5.73e-01
  \\ \hline


   \multirow{8}{*}{\tt Cora\_HA}
   & {\bf Alg.\ref{alg}: 0.67}\newline{}SCFRTR: 4.89\newline{}WYBB: 3.85\newline{}JDCP: 4.07\newline{}RTR: 50.61\newline{}GP-BB: 2.64
   \newline{}PCAL: 2.60\newline{}GPI: 30.50
    &{\bf  1.06e-08}\newline{}7.66e-06\newline{} 4.01e-05\newline{}1.99e-05\newline{}3.24e-07\newline{}1.74e-04  \newline{}7.12e-05\newline{} 3.12e-03
    &{\bf 0}\newline{}1.59e-10\newline{}1.56e-07\newline{}4.07e-09
  \newline{}5.19e-13\newline{}2.65e-06\newline{}2.65e-07
  \newline{}8.64e-04\\ \hline


   \multirow{8}{*}{\tt Cora\_OS} &
    {\bf Alg.\ref{alg}: 1.04}\newline{}SCFRTR: 10.47\newline{}
    WYBB: 10.13\newline{}JDCP: 10.59 \newline{}RTR: 14.49 \newline{}GP-BB: 4.37 \newline{}PCAL: 4.23 \newline{}GPI: 27.02
    & {\bf 8.37e-07}
    \newline{}7.77e-06\newline{}4.16e-05\newline{}3.31e-05\newline{}
    4.52e-06\newline{}1.48e-04\newline{}2.45e-04\newline{}3.40e-03
    &{\bf 0}\newline{}9.89e-10\newline{}5.67e-07\newline{}3.20e-08
  \newline{}1.05e-10\newline{}4.02e-06\newline{}1.85e-05
  \newline{}8.10e-05\\ \hline

   \multirow{8}{*}{\tt Cora\_PL}
    & {\bf Alg.\ref{alg}: 4.64} \newline{} SCFRTR: 17.98\newline{}WYBB: 17.35\newline{}JDCP: 19.18\newline{}
   RTR: 1.73e+02\newline{}GP-BB: 7.47\newline{}PCAL: 7.92
   \newline{}GPI: 64.01
   & {\bf 2.60e-08}\newline{}5.13e-06\newline{} 2.08e-05\newline{}3.73e-05\newline{} 2.95e-07\newline{}1.95e-04 \newline{}2.73e-04\newline{}4.20e-03
   &{\bf 0}\newline{}2.19e-10\newline{}6.44e-08\newline{} 2.04e-07
  \newline{}5.90e-13\newline{}4.39e-06\newline{}1.07e-05
  \newline{}1.81e-05\\ \hline

   \multirow{8}{*}{\tt RCV1\_4Class}
   & {\bf Alg.\ref{alg}: 17.87}\newline{}
   SCFRTR: 1.37e+02\newline{}
   WYBB: 87.703\newline{}
   JDCP: 1.01e+02\newline{}
   RTR: 2.97e+03\newline{}
   GP-BB: 28.58\newline{}PCAL: 20.75\newline{}GPI: 4.48e+02
   &{\bf  8.05e-06}\newline{}9.75e-06\newline{}9.71e-06\newline{}
   9.72e-06\newline{}2.72e-06\newline{}1.15e-04 \newline{}4.04e-04 \newline{}6.84e-05
   & 1.23e-08\newline{}{\bf 0}\newline{}6.73e-08\newline{} 4.45e-08
  \newline{}2.32e-08\newline{}2.42e-06\newline{} 1.70e-05
  \newline{}1.29e-06\\ \hline
    \end{tabular}%
\end{table}%

\subsection{Test on large graph clustering}\label{sec3}
Spectral clustering is a very popular unsupervised machine learning methods \cite{07TSC}.
There are two important problems in spectral clustering  \cite{GCSED}. First, spectral clustering consists of two
successive optimization stages, i.e., spectral embedding and spectral rotation, which may not lead to globally optimal solutions.
Second, for large-scale problems, it is well
known that spectral clustering methods are time-consuming with very high computational complexities. In order to deal with these two challenging problems, a new framework is proposed recently to perform spectral embedding and spectral rotation simultaneously (GCSED) \cite{GCSED}. Unlike the OLSR model, GCSED deals with an $m$-dimensional problem, where $m$ is the number of samples.

Given the database $X=[\bm x_1,\bm x_2,\ldots,\bm x_m]\in\mathbb{R}^{n\times m}$, with $m$ samples and $n$ features drawn from $\ell$ classes. Let
$W$ be a similarity matrix, and $D$ be a diagonal matrix with the diagonal elements being the row sum of $W$. Denote by
$\widehat{W}=D^{-\frac{1}{2}}WD^{-\frac{1}{2}}\in\mathbb{R}^{m\times m}$, and $C=D^{\frac{1}{2}}Y(Y^TDY)^{-\frac{1}{2}}\in \mathbb{R}^{m\times\ell}$,
in each iteration of the GCSED algorithm, one needs to solve the following QMPO problem \cite{GCSED}
\begin{equation}\label{eqn4.5}
\min_{U\in\mathbb{O}^{m\times \ell}} \left\{{\rm tr}\left(U^T(-\widehat{W})U\right)+2\gamma~ {\rm tr}\left(U^T(-{C})\right)   \right\},
\end{equation}
where $Y\in \mathbb{R}^{m\times \ell}$ is the cluster indicator matrix updated in each iteration, and $\gamma$ is a constant trade-off parameter.
Similar to \cite{GCSED}, we make use of the heat kernel weighting to
construct graphs, and calculate edge weights between nodes as
$$
W_{i,j} =\exp\left\{-\frac{\|\bm x_i-\bm x_j\|^2}{2t^2}\right\} ~~{\rm for}~~i,j=1,2,\ldots,m.
$$
In the experiments, we choose $t=0.1$ and $\gamma =0.1,1$, respectively. The data sets used in this example are summarized in Table \ref{Tab2}.



%


The numerical results are reported in Table \ref{fig:2125}, where we run eight algorithms
on the eleven data sets. Some remarks are given. First, we see that all the algorithms are very fast in this example. Second, Algorithm \ref{alg} is the best one in terms of $f_{\rm err}^{(\rm rel)}$ in most of the situations. Third, the proposed algorithm is the best one in terms of KKT error in most of the situations. Indeed, the accuracy of our method can be about five to eight order higher than those of the other methods. Therefore, the proposed block Lanczos method is very promising to large-scale quadratic minimization problems with orthogonality constraints.


\begin{table}[ht]\label{Tab2}
\begin{center}\caption{Summary of test data sets in Example \ref{sec3}.}
\begin{tabular}{l||c|c|c|c}
\hline \hline
 Dataset & \tabincell{c}{\!Number of\\ \!samples ($m$)}& Feature ($n$)&
  \tabincell{c}{\!Number of\\ \!classes ($\ell$)} & Background\\
\hline

\tt AR\tablefootnote{\url{http://www.cad.zju.edu.cn/home/dengcai/Data/TextData.html}.} &1680 &1200 &120 & \multirow{2}{*}{Image}\\
\tt YaleB\tablefootnote{\url{http://cvc.yale.edu/projects/yalefacesB/yalefacesB.html}.} & 2432 & 4069 & 38 \\
\hline
\tt Statlog\tablefootnote{\url{https://archive.ics.uci.edu/ml/machine-learning-databases/statlog/segment/}.} &2310 & 19& 7 &\tabincell{c}{Image\\ segmentation}\\
\hline
\tt madelon\tablefootnote{\url{https://jundongl.github.io/scikit-feature/datasets.html}.}&  2600	&500&	2 &Artificial\\
\hline
\tt TDT2\tablefootnote{The databases {\tt TDT2}, {\tt Reuters} and {\tt 20Newsgroups} are available at \url{http://www.cad.zju.edu.cn/home/dengcai/Data/TextData.html}.}
 & 9394&36771&	30 & Audio\\

\hline
\tt MNIST\tablefootnote{\url{http://yann.lecun.com/exdb/mnist/}.}
 & 70000& 784 &  10& \multirow{3}{*}{\tabincell{c}{ Handwritten\\ text}} \\
\tt USPS\tablefootnote{\url{https://archive.ics.uci.edu/ml/index.php}.} & 9298&256&	10\\
\tt PenDigits\tablefootnote{\url{http://archive.ics.uci.edu/ml}.}
&10992 & 16&	10\\
\hline
\tt Reuters &8293 &18933&65&\multirow{3}{*}{Text}\\
\tt 20Newsgroups&18846 &26214&	20\\
\tt Letters\tablefootnote{\url{https://archive.ics.uci.edu/ml/datasets/Letter+Recognition}.}
&20000	&16&	4\\

\hline\hline
\end{tabular}
\end{center}
\end{table}

\begin{table}[H]
\tiny
  \centering
  \caption{Numerical experiments on the QMPO problem \eqref{eqn4.5} for the data sets in Table \ref{Tab2}, where the best results are in bold.}\label{fig:2125}
    \begin{tabular}{p{8.5em}|p{7.2em}|p{4.34em}|p{4.04em}||p{4.04em}|
    p{4.34em}|p{4.04em}}
    \hline
       &  \multicolumn{3}{c||}{$\gamma=0.1$}   & \multicolumn{3}{c}{$\gamma=1$}  \\
    \hline
       Dataset&  CPU(s) & \tabincell{c}{KKT\\ error }& $f_{\rm err}^{(\rm rel)}$ &CPU(s) & \tabincell{c}{KKT\\ error } & $f_{\rm err}^{\rm (rel)}$\\
    \hline
    \multirow{8}{*}{\tt AR}  & {\bf Alg.\ref{alg}}: 0.56\newline{}SCFRTR: 1.57\newline{}GBB: 0.31\newline{}AFBB: 0.29\newline{}RTR: 0.49\newline{}GP-BB: 0.17\newline{}PCAL: 0.46\newline{}GPI: 1.71 & {\bf 3.11e-09}\newline{}1.29e-06\newline{}5.60e-06\newline{}1.00e-05
    \newline{}1.54e-06\newline{}1.75e-05\newline{}2.93e-05\newline{}1.39e-06&{\bf 0}\newline{}
1.35e-13\newline{}1.43e-12\newline{}7.18e-12\newline{}1.19e-13\newline{}
8.11e-12\newline{}2.55e-11\newline{}2.68e-12&  0.39\newline{}0.87\newline{}0.32\newline{}0.30\newline{} 0.52\newline{}0.14\newline{}0.35\newline{}0.34 &
{\bf 2.71e-13}\newline{}5.50e-06\newline{}1.38e-06\newline{}6.14e-07 \newline{}1.20e-08\newline{}5.39e-06\newline{}1.53e-05\newline{}7.18e-07
  & {\bf 0}\newline{}
9.72e-12\newline{}2.33e-13\newline{}1.06e-13\newline{}
{\bf 0}\newline{}2.76e-12\newline{}5.79e-11\newline{}3.16e-12\\ \hline

\multirow{8}{*}{\tt YaleB} & {\bf Alg.\ref{alg}}: 0.13\newline{}SCFRTR: 0.47\newline{}GBB: 0.14\newline{}AFBB: 0.10\newline{}
RTR: 0.31\newline{}GP-BB: 0.06\newline{}PCAL: 0.10\newline{}GPI: 0.48& {\bf5.82e-09}\newline{}7.14e-07\newline{}3.77e-05\newline{}2.06e-05\newline{}
4.59e-08\newline{}4.79e-05\newline{}6.64e-05\newline{}2.21e-06 &
{\bf 0}\newline{}3.95e-14\newline{}1.19e-10\newline{}1.80e-11\newline{}
1.15e-15\newline{}5.31e-11\newline{}2.82e-10\newline{}3.97e-12& 0.09\newline{}0.22\newline{}0.07\newline{}0.06\newline{}0.10\newline{}
0.03\newline{}0.06\newline{}0.08& {\bf 2.06e-11}\newline{}7.14e-06\newline{}1.46e-06\newline{}4.58e-06\newline{}2.46e-08\newline{}3.73e-06 \newline{}3.71e-06\newline{}9.27e-07&{\bf 0}\newline{}1.42e-11\newline{}
6.34e-13\newline{}3.15e-12\newline{}1.00e-15\newline{}2.00e-12\newline{}
3.15e-12\newline{}1.11e-12\\ \hline

\multirow{8}{*}{\tt Statlog} & {\bf Alg.\ref{alg}}: 0.01\newline{}SCFRTR: 0.02\newline{}GBB: 0.01\newline{}AFBB: 0.01
\newline{}RTR: 0.01\newline{}GP-BB: 0.01\newline{}PCAL: 0.01
\newline{}GPI: 0.05&{\bf 3.14e-09}\newline{}3.51e-07\newline{}2.10e-05
\newline{}3.08e-06\newline{}1.45e-06\newline{}9.85e-06\newline{}5.22e-05\newline{}1.55e-06&
{\bf 0}\newline{}5.13e-15\newline{}1.64e-10\newline{}1.24e-11\newline{}
8.30e-13\newline{}2.46e-11\newline{}2.33e-09\newline{}1.17e-11&
0.01\newline{}0.05\newline{}0.01\newline{}0.01\newline{}0.01
\newline{}0.01\newline{} 0.01\newline{}0.01 & {\bf6.62e-08}\newline{}1.22e-05\newline{}1.09e-06\newline{}
7.18e-06\newline{}1.90e-07\newline{}9.85e-06 \newline{}2.60e-06 \newline{}5.73e-07&{\bf 0}\newline{}3.91e-11\newline{}2.02e-13\newline{}
1.87e-11\newline{}7.45e-15\newline{}4.82e-12\newline{}1.23e-12\newline{}1.61e-13 \\ \hline

\multirow{8}{*}{\tt  madelon} & {\bf Alg.\ref{alg}}: 0.002\newline{}SCFRTR: 0.006\newline{}GBB: 0.007\newline{}AFBB: 0.015\newline{}RTR: 0.007\newline{}GP-BB: 0.004\newline{}PCAL: 0.011\newline{}GPI: 0.027& {\bf 1.72e-13\newline{}2.09e-13}\newline{}8.10e-06\newline{}4.10e-08
\newline{}9.85e-10\newline{}4.20e-08\newline{}9.60e-05 \newline{}9.23e-06 &{\bf 0}\newline{}4.60e-15\newline{}5.48e-12\newline{}1.29e-14\newline{}
1.59e-14\newline{}2.51e-15\newline{}7.68e-10\newline{}1.29e-14
& 0.002\newline{}0.006\newline{}0.003\newline{}0.005\newline{}
0.004\newline{}0.002\newline{}0.002\newline{}0.005 & {\bf 1.58e-14}\newline{} 3.72e-14\newline{}1.43e-08\newline{}1.38e-09
\newline{}4.68e-09\newline{}5.33e-09\newline{}4.57e-05 \newline{}1.07e-06& 8.37e-16\newline{}4.45e-16\newline{}3.76e-15
\newline{}4.60e-15\newline{}{\bf 0}\newline{}5.23e-15\newline{}
6.96e-10\newline{}1.92e-13\\ \hline

\multirow{8}{*}{\tt TDT2} & {\bf Alg.\ref{alg}}: 0.29\newline{}SCFRTR: 0.44\newline{}GBB: 0.28\newline{}AFBB: 0.23\newline{}RTR: 1.79\newline{}GP-BB: 0.16\newline{}PCAL: 0.34\newline{}GPI: 1.29& {\bf 3.43e-14}\newline{}4.79e-06\newline{}7.82e-06\newline{}8.64e-05\newline{}
1.04e-08\newline{}1.15e-04\newline{}1.23e-05 \newline{} 2.44e-06 &{\bf 0}\newline{}1.97e-12\newline{}1.61e-12\newline{}5.12e-13\newline{}6.50e-15\newline{}
1.09e-10\newline{}1.29e-11\newline{}8.40e-12&0.28\newline{} 0.35\newline{}0.20\newline{}0.19\newline{}0.36\newline{}0.10\newline{}0.20\newline{}0.20 & {\bf 2.19e-13}\newline{}2.81e-08\newline{}4.31e-06\newline{}6.03e-07
\newline{}5.64e-11\newline{}2.29e-07\newline{}1.59e-05\newline{}6.87e-07&{\bf 0}\newline{}3.69e-15\newline{}3.87e-12\newline{}4.28e-14\newline{}
2.39e-15\newline{}1.11e-14\newline{}4.29e-11\newline{}9.80e-13\\ \hline

\multirow{8}{*}{\tt MNIST} & {\bf Alg.\ref{alg}}: 0.15\newline{}SCFRTR: 0.64\newline{}GBB: 1.05\newline{}AFBB: 1.56\newline{}RTR: 0.59\newline{}
GP-BB: 0.25\newline{}PCAL: 0.38\newline{} GPI: 4.15& {\bf 1.71e-12\newline{}1.60e-12}\newline{}8.72e-06\newline{}
1.11e-05\newline{}9.57e-12\newline{}1.82e-06\newline{}8.58e-05 \newline{}4.15e-06 &2.43e-14\newline{} 7.33e-14\newline{}6.42e-12\newline{}
1.03e-11\newline{}{\bf 0}\newline{}4.41e-14\newline{}6.13e-10\newline{}
1.16e-11& 0.11\newline{}0.39\newline{}0.44\newline{}0.52\newline{}0.38
\newline{}0.15\newline{}0.23\newline{}0.36& {\bf1.37e-13\newline{}2.96e-13}\newline{}3.98e-06\newline{}
7.71e-06\newline{}9.67e-12\newline{}2.05e-06 \newline{}1.89e-05 \newline{}9.55e-07&9.73e-15\newline{}{\bf 0}\newline{}5.32e-12\newline{}
1.90e-11\newline{}9.36e-15\newline{}3.85e-13\newline{}1.18e-10\newline{}8.54e-13 \\ \hline

\multirow{8}{*}{\tt USPS}  & {\bf Alg.\ref{alg}}: 0.13\newline{}SCFRTR: 0.40\newline{}GBB: 0.08\newline{}AFBB: 0.07\newline{}RTR: 0.11\newline{}GP-BB: 0.05\newline{}PCAL: 0.07\newline{}GPI: 0.28 & {\bf 6.55e-07}\newline{}5.57e-06\newline{}3.61e-05\newline{}4.02e-05\newline{}5.88e-06\newline{}1.48e-05 \newline{}1.77e-04 \newline{}4.50e-06&{\bf 0}\newline{}1.40e-12\newline{}1.11e-10\newline{}1.37e-10\newline{}9.16e-13\newline{}
2.86e-12\newline{}2.66e-09\newline{}1.24e-11&0.08\newline{}0.38\newline{}0.07\newline{}0.05\newline{}
0.05\newline{}0.03\newline{}0.05\newline{}0.06& {\bf 5.96e-08}\newline{}2.81e-05\newline{}2.15e-06\newline{}
7.19e-06\newline{}2.03e-06\newline{}1.73e-05\newline{}6.10e-05\newline{}9.65e-07
&{\bf 0}\newline{}2.45e-10\newline{}1.52e-12\newline{}2.11e-11\newline{}
1.03e-12\newline{}1.24e-10\newline{}6.55e-10\newline{}8.14e-13\\ \hline

\multirow{8}{*}{\tt PenDigits} & {\bf Alg.\ref{alg}}:  0.01\newline{}SCFRTR: 0.04\newline{}GBB: 0.05\newline{}AFBB: 0.06\newline{}RTR: 0.04\newline{}GP-BB: 0.01\newline{}PCAL: 0.02\newline{}GPI: 0.24& {\bf  2.23e-13}\newline{}6.18e-13\newline{}1.26e-06\newline{} 4.78e-08\newline{}2.14e-10\newline{}4.14e-06 \newline{}4.51e-04 \newline{}4.19e-06& {\bf 0}\newline{}1.68e-15\newline{}1.56e-13\newline{}
2.43e-15\newline{}6.55e-15\newline{}3.89e-13\newline{}1.69e-08\newline{}
1.14e-11&0.01\newline{}0.04\newline{}0.08\newline{}0.07\newline{}0.06\newline{}0.02
\newline{}0.03\newline{}0.05& 1.06e-10\newline{}{\bf 8.62e-12}\newline{} 3.21e-06\newline{}2.54e-07\newline{}2.27e-10\newline{}2.79e-06
\newline{}1.55e-04 \newline{}9.71e-07& 2.62e-15\newline{}2.24e-15\newline{}8.68e-13\newline{}1.31e-14\newline{}
{\bf 0}\newline{}7.17e-13\newline{}7.99e-09\newline{}2.62e-15\\ \hline

 \multirow{8}{*}{\tt Reuters} & {\bf Alg.\ref{alg}}: 1.10\newline{}SCFRTR: 5.02\newline{}GBB: 2.78\newline{}AFBB: 1.68\newline{}RTR: 2.36\newline{}GP-BB: 0.66\newline{}PCAL: 1.29\newline{}GPI: 6.35 & {\bf 1.32e-08}\newline{}1.10e-06\newline{}6.62e-05\newline{}1.97e-05
\newline{}2.87e-07\newline{}1.21e-05 \newline{}7.50e-05\newline{}1.88e-06 &{\bf 0} \newline{}6.12e-14 \newline{}3.69e-10 \newline{}9.54e-12 \newline{}7.45e-12 \newline{}3.51e-12 \newline{}4.93e-11 \newline{}7.38e-12&0.40\newline{}2.28\newline{}0.47\newline{}0.47\newline{}0.95\newline{}
0.20\newline{}0.53\newline{}0.48&2.28e-07\newline{}4.16e-06\newline{}5.32e-06\newline{}1.73e-06
\newline{}{\bf 3.15e-09}\newline{}5.00e-05\newline{}5.28e-06 \newline{}5.30e-07& 1.50e-14\newline{}4.62e-12\newline{}
3.30e-12\newline{}3.44e-13\newline{}{\bf 0}\newline{}2.61e-10\newline{}
6.44e-12\newline{}9.72e-13\\ \hline

\multirow{8}{*}{\tt 20Newsgroups} & {\bf Alg.\ref{alg}}: 0.14\newline{}SCFRTR: 0.37\newline{}GBB: 0.82\newline{}AFBB: 0.73\newline{}RTR: 1.36\newline{}GP-BB: 0.22\newline{}PCAL: 0.58\newline{}GPI: 1.63 &  {\bf 9.89e-13}\newline{}2.46e-07\newline{}3.98e-08\newline{}8.64e-05\newline{}
1.22e-06\newline{}4.94e-06\newline{}2.32e-07 \newline{}3.01e-06 &{\bf 0}\newline{}3.44e-15\newline{}1.58e-15\newline{}5.66e-11\newline{}
1.66e-14\newline{}2.03e-12\newline{}5.42e-15\newline{}1.15e-11& 0.18\newline{}0.26\newline{}0.28\newline{}0.31\newline{}0.34\newline{}
0.10\newline{}0.19\newline{}0.28 & {\bf4.70e-13}\newline{}
1.68e-03\newline{}1.08e-05\newline{}4.60e-06\newline{}1.59e-08\newline{}7.73e-06 \newline{}2.64e-05  \newline{}6.87e-07 &{\bf 0}\newline{}
4.71e-07\newline{}1.18e-11\newline{}2.37e-12\newline{}4.10e-15\newline{}
1.80e-11\newline{}2.31e-10\newline{}7.69e-13 \\ \hline

\multirow{8}{*}{\tt Letters} & {\bf Alg.\ref{alg}}: 0.34\newline{}SCFRTR: 2.86\newline{}GBB: 0.82\newline{}AFBB: 0.73
\newline{}RTR: 1.05\newline{}GP-BB: 0.74\newline{}PCAL: 0.58
\newline{}GPI: 2.80 &  {\bf  1.59e-09}\newline{}3.84e-06\newline{}1.41e-05
\newline{}2.63e-05\newline{}8.74e-09\newline{}1.49e-04  \newline{}1.15e-04 \newline{}3.03e-06 &
{\bf 0}\newline{}
1.32e-12\newline{}
1.63e-12\newline{}
5.50e-11\newline{}
4.47e-15\newline{}
2.11e-09\newline{}
1.24e-09\newline{}
8.07e-12 & 0.28\newline{}1.32\newline{}0.36\newline{}0.38\newline{}
    0.52\newline{}0.29\newline{}0.42\newline{}
    0.38 &  {\bf 2.85e-10}\newline{}3.26e-06\newline{}2.41e-05\newline{}
    2.28e-06\newline{}8.87e-06\newline{}2.07e-05\newline{} 5.50e-05\newline{}1.22e-06 &
   {\bf 0}\newline{}
1.60e-12\newline{}
7.58e-11\newline{}
1.32e-12\newline{}
2.30e-11\newline{}
5.37e-11\newline{}
4.17e-10\newline{}
3.02e-12\\ \hline

    \end{tabular}%
  \label{tab:addlabel}%
\end{table}%


\section{Concluding remarks}
In this paper, we propose a block Lanczos method for the large-scale quadratic minimization problems with orthogonality constraints. Convergence analysis on
the optimal value, the optimal solution, the multipliers and the KKT error is given. Theoretical results show that the convergence speed of the new method strictly depends on the distribution of the spectrum of $(I_{\ell}\otimes H)+(\Lambda_*\otimes I_n)$. Specifically, if $(I_{\ell}\otimes H)+(\Lambda_*\otimes I_n)\succ\bf O$, the convergence rate of the solution from the proposed method is comparable to that of conjugate gradient method. Numerical experiments demonstrate that the new algorithm is superior to many state-of-the-art methods for large-scale QMPO in terms of accuracy, KKT error and running time,
especially when $\ell \ll n$.

There are still something deserve further investigation. For instance, as the step $k$ increases, the main workload of Algorithm \ref{alg} is to solve \eqref{3}. The computational overhead will be prohibitive if $k$ is large, and we have to restrict the value of $k$, and efficient restarting techniques \cite{2,Y.S,Stewart} are required for our block Krylov subspace method.
On the other hand, we assume that the matrix $(I_{\ell}\otimes H)+(\Lambda_*\otimes I_n)$ is nonsingular in the convergence analysis.
An interesting topic is to weaken this assumption for the analysis.

%

%


\section*{Acknowledgement}
We are grateful to Prof. Leihong Zhang and Prof. Chungen Shen for providing us some codes and databases used in numerical experiments. Meanwhile, we would like to thank Dr. Yongyan Guo for helpful discussions.

\end{document}